\documentclass[12pt]{article}
\usepackage{amssymb}
\usepackage{latexsym}
\usepackage{amscd}

\newtheorem{definition}{Definition}[section]
\newtheorem{theorem}[definition]{Theorem}
\newtheorem{lemma}[definition]{Lemma}

\newtheorem{note}[definition]{Note}

\newtheorem{notation}[definition]{Notation}

\begin{document}

\title{\bf Tridiagonal pairs of   \\
shape $(1, 2, 1)$
}
\author{
Melvin A. Vidar{\footnote{
Math and Statistics Department,
College of Arts and Sciences,
University of the East-Manila, Philippines.
}}}

\date{}

\maketitle
\begin{abstract}
Let $\mathbb F$ denote a field and let $V$ denote a vector space over
$\mathbb F$ with finite positive dimension.
We consider a pair of linear transformations
$A:V\to V$ and $A^*:V\to V$ that satisfies the following conditions:
(i) each of $A,A^*$ is diagonalizable; (ii) there exists an ordering
$\lbrace V_i \rbrace_{i=0}^d$ of the eigenspaces of $A$ such that
$A^* V_i \subseteq V_{i-1}+V_i+V_{i+1}$ for $0 \leq i \leq d$,
where $V_{-1} = 0$ and $V_{d+1} = 0$;
(iii) there exists an ordering
$\lbrace V^*_i \rbrace_{i=0}^{\delta}$ of the eigenspaces of $A^*$ such that
$AV^*_i \subseteq V^*_{i-1}+V^*_i+V^*_{i+1}$ for $0 \leq i \leq \delta $,
where $V^*_{-1} = 0$ and $V^*_{\delta+1} = 0$;
(iv) there is no subspace $W$ of $V$ such that $AW\subseteq W$, 
 $A^*W\subseteq W$, $W \neq 0, W \neq V$. We call such a pair a {\it tridiagonal pair}
on $V$. It is known that $d = \delta$ and that 
for $0 \leq i \leq d$
the dimensions of $V_i, V_{d-i}, V^*_i, V^*_{d-i}$ coincide;
we denote this common value by $\rho_i$. The sequence
$\lbrace \rho_i\rbrace_{i=0}^d$ is called the {\it shape} of
the pair. In this paper we assume the shape is $(1,2,1)$ and
obtain the following results. We describe six bases for
$V$; one diagonalizes $A$, another diagonalizes $A^*$, and the
other four underlie the split decompositions for $A,A^*$.
We give the action of $A$ and $A^*$ on each basis. For each
ordered pair of bases among the six, we give the
transition matrix. At the end we classify the tridiagonal
pairs of shape $(1,2,1)$ in terms of a sequence of scalars
called the parameter array.

\end{abstract}

\section{Introduction}
\noindent Throughout this paper $\mathbb F$ will denote a field and $V$ will be a vector space over $\mathbb F$ with finite positive dimension.\\

\medskip
\noindent By a {\it linear transformation on $V$} we mean an $\mathbb F$-linear map from $V$ to $V$. Let $A$ denote a linear transformation on $V$. By an {\it eigenspace} of $A$ we mean a nonzero subspace of $V$ of the form
\begin{equation}
\lbrace v \in V \;\vert \;Av = \theta v\rbrace,
\label{eq:defeigspace}
\end{equation}
where $\theta \in \mathbb F$. We say $A$ is {\it diagonalizable} on $V$ whenever $V$ is spanned by the eigenspaces of $A$.

\begin{definition} \rm \cite[Definition 1.1]{pandq}
\label{def:pnq} By a {\it tridiagonal pair} (or TD pair) on $V$, we mean an
ordered pair $(A, A^*)$, where $A$ and $A^*$ are linear
transformations on $V$ that satisfy the following four conditions.

\begin{enumerate}
\item $A$ and $A^*$ are both diagonalizable on $V$.
\item There exists an ordering $\lbrace V_i\rbrace_{i=0}^d$ of the
eigenspaces of $A$ such that
\begin{eqnarray}
\label{eq:TD1} A^*V_i \subseteq V_{i-1} + V_i + V_{i+1},\qquad (0
\leq i\leq d),
\end{eqnarray}
where $V_{-1}= 0, V_{d+1} = 0$.
\item There exists an ordering $\lbrace V^*_i\rbrace_{i=0}^\delta$ of the
eigenspaces of $A^*$ such that
\begin{eqnarray}
\label{eq:TD2} AV^*_i \subseteq V^*_{i-1} + V^*_i +
V^*_{i+1},\qquad (0 \leq i\leq \delta),
\end{eqnarray}
where $V^*_{-1}= 0$, $V^*_{\delta+1} = 0.$
\item There is no subspace $W$ of $V$ such that $AW \subseteq
W$, $A^*W \subseteq W$, $W \neq 0$, $W \neq V$.
\end{enumerate}
\end{definition}

\begin{note}\rm
According to common notational convention $A^*$ denotes the conjugate-transpose of $A$.
We emphasize that we are not using this convention. In a TD
pair $(A, A^*)$ the linear transformations $A$ and $A^*$ are
arbitrary subject to (i)--(iv) above.
\end{note}

\section{TD systems}
When working with a TD pair, it is often convenient to consider a
closely related but somewhat more abstract object called a TD
system. To define it, we recall a few concepts from linear
algebra. Let End($V$) denote the $\mathbb F$-algebra
consisting of all linear transformations on $V$ and let $A$ denote a
diagonalizable element in End($V$). Let $\lbrace \theta_i\rbrace_{i=0}^d$
denote an ordering of the eigenvalues of $A$, and put
\begin{eqnarray}
\label{eq:Eidef}
E_i = \prod_{{0 \leq  j \leq d}\atop {j\not=i}} {{A-\theta_j
I}\over {\theta_i-\theta_j}} \label{eq:idem3}
\end{eqnarray}

\noindent for $0 \leq i \leq d.$ By elementary linear algebra,

\begin{eqnarray}
\label{eq:idem4} AE_i = E_iA &=& \theta_iE_i \qquad (0 \leq i \leq d),\\
\label{eq:idem5} E_iE_j &=& \delta_{ij}E_i \qquad (0 \leq i, j \leq d),\\
\label{eq:idem6} \sum_{i=0}^d E_i &=& I.
\end{eqnarray}

\noindent From this, one finds $\lbrace E_i\rbrace_{i=0}^d$ is a basis for the
subalgebra of End($V$) generated by $A$. We refer to $E_i$ as the
{\it primitive idempotent} of $A$ associated with $\theta_i$. It
is helpful to think of these primitive idempotents as follows. From
(\ref{eq:idem5}), (\ref{eq:idem6}) one readily finds
\\
\begin{eqnarray*}
\label{eq:eq7} V = E_0V + E_1V +\cdots+ E_dV \qquad \qquad
(\mbox{direct sum}).
\end{eqnarray*}
\\
\noindent For $0\leq i \leq d$, $E_iV$ is the eigenspace of $A$ in $V$ associated with the eigenvalue
$\theta_i$ and $E_i$ acts  on $V$ as the projection onto this
eigenspace.

\medskip

\begin{definition} \rm \cite[Definition 2.1]{pandq}
\label{def:TDsys} By a {\it tridiagonal system} (or TD system) on $V$, we mean
a sequence
\begin{eqnarray}
\label{eq:eq8} \Phi:=(A; {\lbrace{E_i}\rbrace}_{i=0}^d; A^*;{\lbrace{E_i^*}
\rbrace}_{i=0}^\delta)
\end{eqnarray}
\noindent that satisfies (i)--(vi) below.

\begin{enumerate}
\item $A$ and $A^*$ are both diagonalizable linear transformations on $V$.
\item $\lbrace E_i\rbrace_{i=0}^d$ is an ordering of the primitive
idempotents of $A$.
\item $\lbrace E^*_i\rbrace_{i=0}^\delta$ is an ordering of the
primitive idempotents of $A^*$.
\item $E_iA^*E_j = 0$ if $ \mid i-j \mid > 1,\qquad (0 \leq
i,j\leq d)$.
\item $E^*_iAE^*_j = 0$ if $ \mid i-j \mid > 1,\qquad (0 \leq
i,j\leq \delta)$.
\item There is no subspace $W$ of $V$ such that $AW \subseteq
W$, $A^*W \subseteq W$, $W \neq 0$, $W \neq V$.
\end{enumerate}
\end{definition}
\medskip

\begin{lemma} {\rm \cite[Lemma 2.2]{pandq}}
\label{lem:TDPhi} Let $(A, A^*)$ denote a TD pair on $V$. Let ${\lbrace{V_i}
\rbrace}_{i=0}^d$ denote an ordering of the eigenspaces of $A$ satisfying
(\ref{eq:TD1}) and for $0 \leq i \leq d$ let $E_i$ denote the
primitive idempotent of $A$ associated with $V_i$. Let ${\lbrace{V^*_i}
\rbrace}_{i=0}^\delta$ denote an ordering of the eigenspaces of $A^*$ satisfying
(\ref{eq:TD2}) and for $0 \leq i \leq \delta$ let $E^*_i$ denote the the primitive idempotent of $A^*$ associated with $V^*_i$. Then $(A; \lbrace{E_i}
\rbrace_{i=0}^d; A^*;\lbrace{E_i^*}
\rbrace_{i=0}^\delta)$
is a TD system on $V$.
\end{lemma}

\medskip

\noindent Let $\Phi:=(A; {\lbrace{E_i}\rbrace}_{i=0}^d; A^*;{\lbrace{E_i^*}
\rbrace}_{i=0}^\delta)$ denote a TD system on $V$. By \cite[Lemma 2.3]{pandq} the pair $(A, A^*)$ is a TD pair on $V$; we say this pair is {\it associated with} $\Phi$.\\

\medskip
\noindent Referring to Definition \ref{def:pnq} and Definition \ref{def:TDsys}, it turns out that  $d = \delta$ \cite[Lemma 4.5]{pandq}; we call this common value the {\it diameter}.

\section{The relatives of a TD system}

A given TD system  can be modified in  several ways to get a new TD system. For instance,
let $\Phi$ denote the TD system from Definition \ref{def:TDsys}. Then each of
\begin{eqnarray}
\label{Phistar} \Phi^*&:=&(A^*;\lbrace{E_i^*}\rbrace_{i=0}^d; A; {\lbrace{E_i}\rbrace}_{i=0}^d),\\
\label{Phidown} \Phi^\downarrow&:=&(A; \lbrace{E_i}\rbrace_{i=0}^d; A^*;\lbrace E^*_{d-i}\rbrace_{i=0}^d),\\
\label {Phiddowm} \Phi^\Downarrow&:=&(A; \lbrace E_{d-i}\rbrace_{i=0}^d; A^*;\lbrace{E_i^*}\rbrace_{i=0}^d)
\end{eqnarray}

\noindent is a TD system on $V$. Viewing $*, \downarrow, \Downarrow$ as permutations on the set of all TD systems,
\begin{eqnarray}
&&\qquad \qquad \qquad  *^2 \;=\;
\downarrow^2\;= \;
\Downarrow^2 \;=\;1,
\qquad \quad
\label{eq:deightrelationsAS99}
\\
&&\Downarrow *\;
=\;
* \downarrow,\qquad \qquad
\downarrow *\;
=\;
* \Downarrow,\qquad \qquad
\downarrow \Downarrow \; = \;
\Downarrow \downarrow.
\qquad \quad
\label{eq:deightrelationsBS99}
\end{eqnarray}
The group generated by symbols
$*, \downarrow, \Downarrow $ subject to the relations
(\ref{eq:deightrelationsAS99}),
(\ref{eq:deightrelationsBS99})
is the dihedral group $D_4$.
We recall $D_4$ is the group of symmetries of a square,
and has 8 elements.
Apparently $*, \downarrow, \Downarrow $ induce an action of
 $D_4$ on the set of all TD systems.
Two TD systems will be called {\it relatives} whenever they
are in the same orbit of this $D_4$ action.
The relatives of $\Phi$ are as follows:
\medskip

\centerline{
\begin{tabular}[t]{c|c}
        name &relative \\ \hline
        $\Phi$ & $(A;{\lbrace{E_i}\rbrace}_{i=0}^d;A^*;{\lbrace{E^*_i}\rbrace}_{i=0}^d)$   \\
        $\Phi^\downarrow$ & $(A;{\lbrace{E_i}\rbrace}_{i=0}^d;A^*;{\lbrace E^*_{d-i}\rbrace}_{i=0}^d)$   \\
        $\Phi^\Downarrow$ & $(A;{\lbrace E_{d-i}\rbrace}_{i=0}^d;A^*;{\lbrace{E^*_i}\rbrace}_{i=0}^d)$   \\
        $\Phi^{\downarrow \Downarrow}$ &$(A;{\lbrace E_{d-i}\rbrace}_{i=0}^d;A^*;{\lbrace E^*_{d-i}\rbrace}_{i=0}^d)$   \\
    $\Phi^*$ & $(A^*;{\lbrace{E^*_i}\rbrace}_{i=0}^d;A;{\lbrace{E_i}\rbrace}_{i=0}^d)$   \\
        $\Phi^{\downarrow *}$ & $(A^*;{\lbrace E^*_{d-i}\rbrace}_{i=0}^d;A;{\lbrace{E_i}\rbrace}_{i=0}^d)$ \\
        $\Phi^{\Downarrow *}$ & $(A^*;{\lbrace{E^*_i}\rbrace}_{i=0}^d;A;{\lbrace E_{d-i}\rbrace}_{i=0}^d)$ \\
    $\Phi^{\downarrow \Downarrow *}$ & $(A^*;{\lbrace E^*_{d-i}\rbrace}_{i=0}^d;A;{\lbrace E_{d-i}\rbrace}_{i=0}^d)$
\end{tabular}}

\bigskip

\noindent We now introduce two sequences of parameters that we will use to describe a given TD system.\\

\begin{definition}\rm \cite[Definition 3.1]{pandq}
\label{def:eigenseq} Let $\Phi:=(A; {\lbrace{E_i}\rbrace}_{i=0}^d; A^*;{\lbrace{E_i^*}
\rbrace}_{i=0}^d)$ denote a TD system.
For $0 \leq i \leq d$ let $\theta_i$ denote the eigenvalue of $A$ associated with $E_i$.  We refer to
${\lbrace{\theta_i}\rbrace}_{i=0}^d$ as the {\it eigenvalue sequence of $\Phi$}. For $0 \leq i \leq d$ let $\theta^*_i$ denote the eigenvalue of $A^*$ associated with $E^*_i$. We refer to ${\lbrace{\theta^*_i}\rbrace}_{i=0}^d$ as the {\it dual eigenvalue sequence of $\Phi$}. We remark that $\theta_0, \theta_1, \ldots,\theta_d$ are mutually distinct and  $\theta^*_0, \theta^*_1, \ldots, \theta^*_d$ are mutually distinct.
\end{definition}

\medskip
\section{The split decomposition}

\begin{definition}\rm By a {\it decomposition of $V$ of length $d$} we mean a sequence
${\lbrace{V_i}\rbrace}_{i=0}^d$ of nonzero subspaces of $V$ such that
\begin{eqnarray*}
V = V_0 + V_1 +\cdots+ V_d \qquad \qquad
(\mbox{direct sum}).
\end{eqnarray*}
\end{definition}

\noindent We do not assume each of $V_0, V_1,\ldots,V_d$ has dimension 1. For $0 \leq i \leq d$ we call $V_i$ the {\it ith component} of the decomposition. For notational convenience we let $V_{-1}$ = 0 and $V_{d + 1}$ = 0. \\

\noindent We will refer to the following setup.

\medskip
\begin{notation} \rm \label{note:decomp} In this section we let
\begin{eqnarray*}
\Phi:=(A; {\lbrace{E_i}\rbrace}_{i=0}^d; A^*;{\lbrace{E_i^*}\rbrace}_{i=0}^d)
\end{eqnarray*}
be a TD system on $V$ with eigenvalue sequence ${\lbrace{\theta_i}\rbrace}_{i=0}^d$ and dual eigenvalue sequence ${\lbrace{\theta^*_i}\rbrace}_{i=0}^d$.
\end{notation}

\medskip
\noindent With reference to Notation \ref{note:decomp}, we are about to define six decompositions of $V$.
In order to keep track of these decompositions we will give each of them a name. Our naming scheme is as follows. Let $\Omega$ denote the set consisting of the four symbols $0, D, 0^*, D^*$. Each of the six decompositions will get a name $[u]$ where $u$ is a two element subset of $\Omega$.

\bigskip
\begin{lemma} {\rm \cite[Lemma 4.2]{quantumaff}} \label{lem:sixdecomp}
With reference to Notation \ref{note:decomp}, for each of the six rows in the table below and for $0 \leq i \leq d$ let $U_i$ denote the ith component described in that row. Then the sequence $\lbrace U_i\rbrace_{i=0}^d$ is a decomposition of V.
\medskip

\centerline{
\begin{tabular}[t]{c|c}
        decomposition &$ith$ component \\ \hline
    $[0^*D]$ & $(E^*_0V + \cdots+E^*_iV) \cap (E_iV + \cdots + E_dV)$ \\
        $[0^*0]$ & $(E^*_0V + \cdots+E^*_iV) \cap (E_0V + \cdots + E_{d-i}V)$ \\
        $[D^*0]$ & $(E^*_dV + \cdots+E^*_{d-i}V) \cap (E_0V + \cdots + E_{d-i}V)$ \\
        $[D^*D]$ & $(E^*_dV + \cdots+E^*_{d-i}V) \cap (E_iV + \cdots + E_dV)$ \\
    $[0D]$ & $E_iV$ \\
    $[0^*D^*]$ & $E^*_iV$ \\
\end{tabular}}
\end{lemma}

\noindent Referring to the table in Lemma \ref{lem:sixdecomp}, we call the decompositions corresponding to the first four rows the {\it split decompositions} of $V$. We observe that the last two rows give the eigenspace decompositions of $A$ and $A^*$.

\medskip
\begin{lemma} {\rm \cite [Lemma 5.1]{quantumaff}} \label{lem:actionsplit}
With reference to Notation \ref{note:decomp}, let ${\lbrace{U_i}\rbrace}_{i=0}^d$ denote any one of the six decompositions of $V$ given in Lemma \ref{lem:sixdecomp}. Then for $0 \leq i \leq d$ the action of $A$ and $A^*$ on $U_i$ is described as follows.

\medskip
\centerline{
\begin{tabular}[t]{c|c|c}
        decomposition &action of $A$ on $U_i$&action of $A^*$ on $U_i$ \\ \hline
        $[0^*D]$ & $(A - \theta_iI)U_i \subseteq U_{i+1}$& $(A^* - \theta^*_iI)U_i \subseteq U_{i-1}$ \\
        $[0^*0]$ & $(A - \theta_{d-i}I)U_i \subseteq U_{i+1}$& $(A^* - \theta^*_iI)U_i \subseteq U_{i-1}$ \\
        $[D^*0]$ & $(A - \theta_{d-i}I)U_i \subseteq U_{i+1}$& $(A^* - \theta^*_{d-i}I)U_i \subseteq U_{i-1}$ \\
        $[D^*D]$ & $(A - \theta_iI)U_i \subseteq U_{i+1}$& $(A^* - \theta^*_{d-i}I)U_i \subseteq U_{i-1}$ \\
        $[0D]$ & $(A - \theta_iI)U_i = 0$ & $A^*U_i \subseteq U_{i-1} + U_i + U_{i+1}$ \\
        $[0^*D^*]$ & $AU_i \subseteq U_{i-1} + U_i + U_{i+1}$& $(A^* - \theta^*_iI)U_i = 0$ \\
    \end{tabular}}
\end{lemma}

\medskip \begin{lemma} {\rm \cite[Lemma 4.4]{quantumaff}} \label{lem:shape}
With reference to Notation \ref{note:decomp}, let ${\lbrace{U_i}\rbrace}_{i=0}^d$ denote any one of the six decompositions of $V$ given in Lemma \ref{lem:sixdecomp}. For $0 \leq i\leq d$ let $\rho_i$ denote the dimension of $U_i$. Then the sequence ${\lbrace{\rho_i}\rbrace}_{i=0}^d$ is independent of the decomposition. Moreover, this sequence is unimodal and symmetric; that is $\rho_i = \rho_{d-i}$ for $0 \leq i \leq d$ and $\rho_{i-1} \leq \rho_i$ for $1 \leq i \leq d/2$.
\end{lemma}

\noindent Referring to Lemma \ref{lem:shape}, we call the sequence ${\lbrace{\rho_i}\rbrace}_{i=0}^d$ the {\it shape} of the TD system $\Phi$. A TD system of shape $(1, 1, \ldots, 1)$ is the same as a {\it Leonard system} \cite{LPPaper}.

\section{Some parameters}

For the rest of the paper we are going to consider a TD system of diameter 2 and shape $(1, 2, 1)$. We refer to the following setup.

\begin{notation} \rm \label{note:notephi} Fix a TD system $\Phi = (A; {\lbrace{E_i}\rbrace}_{i=0}^2; A^*;{\lbrace{E^*_i}\rbrace}_{i=0}^2)$ on $V$ with eigenvalue sequence ${\lbrace{\theta_i}\rbrace}_{i=0}^2$, dual eigenvalue sequence ${\lbrace{\theta^*_i}\rbrace}_{i=0}^2$, and shape $(1, 2, 1)$.\\

\noindent Setting $d = 2$ in (\ref{eq:idem3}) we obtain the following elements in End($V$):
\begin{eqnarray*}
E_0 = \frac{(A - \theta_1I)(A - \theta_2I)}{(\theta_0 - \theta_1)(\theta_0 - \theta_2)}, \qquad \qquad
E^*_0 = \frac{(A^* - \theta^*_1I)(A^* - \theta^*_2I)}{(\theta^*_0 - \theta^*_1)(\theta^*_0 - \theta^*_2)}
\end{eqnarray*}

\begin{eqnarray*}
E_1 = \frac{(A - \theta_0I)(A - \theta_2I)}{(\theta_1 - \theta_0)(\theta_1 - \theta_2)}, \qquad \qquad
E^*_1 = \frac{(A^* - \theta^*_0I)(A^* - \theta^*_2I)}{(\theta^*_1 - \theta^*_0)(\theta^*_1 - \theta^*_2)}
\end{eqnarray*}

\begin{eqnarray*}
E_2 = \frac{(A - \theta_0I)(A - \theta_1I)}{(\theta_2 - \theta_0)(\theta_2 - \theta_1)}, \qquad \qquad
E^*_2 = \frac{(A^* - \theta^*_0I)(A^* - \theta^*_1I)}{(\theta^*_2 - \theta^*_0)(\theta^*_2 - \theta^*_1)}.
\end{eqnarray*}
\end{notation}

\bigskip
\noindent Referring to Notation \ref{note:notephi}, in order to describe $\Phi$ we will need two parameters $\varphi$ and $\phi$. We now introduce these parameters. \\

\medskip \begin{lemma}
\label{lem:E0stareigen}
With reference to Notation \ref{note:notephi}, there exist nonzero scalars $\varphi, \phi \in \mathbb F$ such that (i), (ii) hold below.
\begin{enumerate}
\item $E^*_0V$ is an eigenspace for $(A^* -\theta^*_1I)(A^* -\theta^*_2I)(A - \theta_1I)(A - \theta_0I)$ and the corresponding eigenvalue is $\varphi$.
\item $E^*_0V$ is an eigenspace for $(A^* -\theta^*_1I)(A^* -\theta^*_2I)(A - \theta_1I)(A - \theta_2I)$ and the corresponding eigenvalue is $\phi$.
\end{enumerate}

\noindent \rm {\bf Proof}. (i) Referring to row $[0^*D]$ of Lemma \ref{lem:actionsplit} we have
$(A - \theta_1I)(A - \theta_0I)U_0 \subseteq U_2$ and $(A^* - \theta^*_1I)(A^* - \theta^*_2I)U_2 \subseteq U_0.$ In both cases equality is obtained by \cite[Lemma 6.5]{pandq}. By these comments $U_0$ is an eigenspace for $(A^* - \theta^*_1I)(A^* - \theta^*_2I)(A - \theta_1I)(A - \theta_0I)$ and the corresponding eigenvalue is a nonzero scalar in $\mathbb F$. We denote this eigenvalue by $\varphi$. By row $[0^*D]$ of Lemma \ref{lem:sixdecomp} we have $U_0 = E^*_0V$ and the result follows.\\

\noindent (ii) Referring to row $[0^*0]$ of Lemma \ref{lem:actionsplit} we have $(A - \theta_1I)(A - \theta_2I)U_0 \subseteq U_2$ and $(A^* - \theta^*_1I)(A^* - \theta^*_2I)U_2 \subseteq U_0$. In both cases equality is attained by \cite[Lemma 6.5]{pandq}. By these comments $U_0$ is an eigenspace for $(A^* - \theta^*_1I)(A^* - \theta^*_2I)(A - \theta_1I)(A - \theta_2I)$ and the corresponding eigenvalue is a nonzero scalar in $\mathbb F$. We denote this eigenvalue by $\phi$. By row $[0^*0]$ of Lemma \ref{lem:sixdecomp} we have $U_0 = E^*_0V$ and the result follows. \hfill $\Box$\\
\end{lemma}

\medskip \begin{lemma}
\label{lem:E0eigen}
With reference to Notation \ref{note:notephi} and Lemma \ref{lem:E0stareigen}, the following (i), (ii) hold.
\begin{enumerate}
\item $E_0V$ is an eigenspace for $(A - \theta_1I)(A - \theta_2I)(A^* -\theta^*_1I)(A^* -\theta^*_2I)$ and the corresponding eigenvalue is $\phi$.
\item $E_0V$ is an eigenspace for $(A - \theta_1I)(A - \theta_2I)(A^* -\theta^*_1I)(A^* -\theta^*_0I)$ and the corresponding eigenvalue is $\varphi$.
\end{enumerate}

\noindent \rm {\bf Proof}. (i) Let $0 \neq u_0 \in E_0V$. By \cite[Lemma 6.5]{pandq} there exists $0 \neq u^*_0 \in E^*_0V$ such that $(A - \theta_1I)(A - \theta_2I)u^*_0 = u_0$. Combining this fact with Lemma \ref{lem:E0stareigen}(ii) we get

\begin {eqnarray*}
&& (A - \theta_1I)(A - \theta_2I)(A^* -\theta^*_1I)(A^* -\theta^*_2I)u_0\\
&& \qquad \qquad = \quad (A - \theta_1I)(A - \theta_2I)(A^* -\theta^*_1I)(A^* -\theta^*_2I)(A - \theta_1I)(A - \theta_2I)u^*_0\\
&& \qquad \qquad = \quad \phi (A - \theta_1I)(A - \theta_2I)u^*_0\\
&& \qquad \qquad = \quad \phi u_0.
\end{eqnarray*}
\noindent The result follows.\\
\noindent (ii) Pick $0 \neq u_0 \in E_0V$. By \cite[Lemma 6.5]{pandq} there exists $0 \neq u^*_0 \in E^*_0V$ such that $(A - \theta_1I)(A - \theta_2I)u^*_0 = u_0$. By Lemma \ref{lem:E0stareigen}(i),
\begin{eqnarray*}
(A^* - \theta^*_1I)(A^* - \theta^*_2I)(A - \theta_1I)(A - \theta_0I)u^*_0
&=& \varphi u^*_0.
\end{eqnarray*}

\noindent By row $[0^*D^*]$ of Lemma \ref{lem:actionsplit},
\begin{eqnarray*}
(A^* - \theta^*_1I)(A^* - \theta^*_0I)(A - \theta_1I)u^*_0 &=& 0.\\
\end{eqnarray*}

\noindent By rows $[0^*D]$ and $[0D]$  of Lemma \ref{lem:actionsplit},
\begin{eqnarray*}
(A - \theta_1I)(A - \theta_2I)(A^* - \theta^*_1I)(A - \theta_1I)(A - \theta_0I)u^*_0 &=& 0.\\
\end{eqnarray*}

\medskip
\noindent Observe that
\begin{eqnarray*}
&& (A^* - \theta^*_0I)(A - \theta_1I)(A - \theta_2I)\\
&& \qquad \qquad = \quad (A^* - \theta^*_2I)(A - \theta_1I)(A - \theta_0I)\\
&& \qquad \qquad \quad \; + \quad (\theta^*_2 - \theta^*_0)(A - \theta_1I)(A - \theta_0I)\\
&& \qquad \qquad  \quad \; + \quad (\theta_0 - \theta_2)(A^* - \theta^*_0I)(A - \theta_1I).
\end{eqnarray*}
\noindent Using the above lines we argue
\begin{eqnarray*}
&&(A - \theta_1I)(A - \theta_2I)(A^* -\theta^*_1I)(A^* -\theta^*_0I)u_0\\
&&\qquad \qquad = \quad (A - \theta_1I)(A - \theta_2I)(A^* - \theta^*_1I)(A^* - \theta^*_0I)(A - \theta_1I)(A - \theta_2I)u^*_0 \\
&&\qquad \qquad = \quad(A - \theta_1I)(A - \theta_2I)(A^* - \theta^*_1I)(A^* - \theta^*_2I)(A - \theta_1I)(A - \theta_0I)u^*_0 \\
&& \qquad \qquad \quad \; + \quad (\theta^*_2 - \theta^*_0) (A - \theta_1I)(A - \theta_2I)(A^* - \theta^*_1I)(A - \theta_1I)(A - \theta_0I)u^*_0 \\
&&\qquad \qquad \quad \; + \quad (\theta_0 - \theta_2)(A - \theta_1I)(A - \theta_2I)(A^* - \theta^*_1I)(A^* - \theta^*_0I)(A - \theta_1I)u^*_0 \\
&&\qquad \qquad = \quad \varphi (A - \theta_1I)(A - \theta_2I) u^*_0 \quad + \quad 0 \quad + \quad0\\
&&\qquad \qquad = \quad \varphi u_0
\end{eqnarray*}

\noindent and the result follows. \hfill $\Box$
\end{lemma}

\medskip \begin{definition} \label{def:pararray}\rm With reference to Notation \ref{note:notephi}, by the {\it first split eigenvalue} (respectively {\it second split eigenvalue}) for $\Phi$ we mean the scalar $\varphi$ (respectively $\phi$) in Lemma \ref{lem:E0stareigen}. By the {\it parameter array} of $\Phi$ we mean the sequence $({\lbrace{\theta_i}\rbrace}_{i=0}^2; \lbrace{\theta_i^*}\rbrace_{i=0}^2; \varphi; \phi)$.
\end{definition}

\medskip \begin{lemma}
\label{lem:relpararray} With reference to Definition \ref{def:pararray} the following (i)--(iii) hold.
\begin{enumerate}
\item The parameter array of $\Phi^*$ is $({\lbrace{\theta^*_i}\rbrace}_{i=0}^2; \lbrace{\theta_i}\rbrace_{i=0}^2; \varphi; \phi)$.
\item The parameter array of $\Phi^\downarrow$ is $({\lbrace{\theta_i}\rbrace}_{i=0}^2; \lbrace{\theta_{2-i}^*}\rbrace_{i=0}^2; \phi; \varphi)$.
\item The parameter array of $\Phi^\Downarrow$ is $({\lbrace{\theta_{2-i}}\rbrace}_{i=0}^2; \lbrace{\theta_i^*}\rbrace_{i=0}^2; \phi; \varphi)$.
\end{enumerate}

\noindent \rm {\bf Proof}. (i) Apply Lemma \ref{lem:E0stareigen} to $\Phi^*$ and use Lemma \ref{lem:E0eigen}.\\
\noindent (ii) Apply Lemma \ref{lem:E0eigen} to $\Phi^\downarrow$.\\
\noindent (iii) Apply Lemma \ref{lem:E0stareigen} to $\Phi^\Downarrow$. \hfill $\Box$
\end{lemma}

\medskip \begin{lemma}
\label{lem:E2stareigen}
With reference to Notation \ref{note:notephi} and Lemma \ref{lem:E0stareigen}, the following (i), (ii) hold.
\begin{enumerate}
\item $E^*_2V$ is an eigenspace for $(A^* -\theta^*_1I)(A^* -\theta^*_0I)(A - \theta_1I)(A - \theta_2I)$ and the corresponding eigenvalue is $\varphi$.
\item $E^*_2V$ is an eigenspace for $(A^* -\theta^*_1I)(A^* -\theta^*_0I)(A - \theta_1I)(A - \theta_0I)$ and the corresponding eigenvalue is $\phi$.
\end{enumerate}

\noindent \rm {\bf Proof}. Apply Lemma  \ref{lem:E0stareigen}(i), (ii) to $\Phi^\downarrow$ and use Lemma \ref{lem:relpararray}(ii). \hfill $\Box$\\
\end{lemma}

\medskip \begin{lemma}
\label{lem:E2eigen}
With reference to Notation \ref{note:notephi} and Lemma \ref{lem:E0stareigen}, the following (i), (ii) hold.
\begin{enumerate}
\item $E_2V$ is an eigenspace for $(A - \theta_1I)(A - \theta_0I)(A^* - \theta^*_1I)(A^* - \theta^*_2I)$ and the corresponding eigenvalue is $\varphi$.
\item $E_2V$ is an eigenspace for $(A - \theta_1I)(A - \theta_0I)(A^* - \theta^*_1I)(A^* - \theta^*_0I)$  and the corresponding eigenvalue is $\phi$.
\end{enumerate}

\noindent \rm {\bf Proof}. Apply Lemma  \ref{lem:E0eigen}(i), (ii) to $\Phi^\Downarrow$ and use Lemma \ref{lem:relpararray}(iii). \hfill $\Box$
\end{lemma}

\medskip \begin{lemma} \label{lem:releigenseq}
With reference to Notation \ref{note:notephi},
\begin{enumerate}
\item $E^*_0V$ is an eigenspace for $(A^* - \theta^*_1I)(A - \theta_0I)$ and the corresponding eigenvalue is
\begin{eqnarray}
\label{eq:lambda}
\varphi_1:= \frac{\phi - \varphi}{(\theta_0 - \theta_2)(\theta^*_0 - \theta^*_2)} -  (\theta_0 - \theta_1)(\theta^*_0 - \theta^*_1).
\end{eqnarray}

\item $E^*_0V$ is an eigenspace for $(A^* - \theta^*_1I)(A - \theta_2I)$ and the corresponding eigenvalue is
\begin{eqnarray}
\label{eq:lambda2}
\phi_1:= \frac{\varphi - \phi}{(\theta_2 - \theta_0)(\theta^*_0 - \theta^*_2)} -  (\theta_2 - \theta_1)(\theta^*_0 - \theta^*_1).
\end{eqnarray}

\item $E^*_2V$ is an eigenspace for $(A^* - \theta^*_1I)(A - \theta_0I)$ and the corresponding eigenvalue is
\begin{eqnarray}
\label{eq:mu2}
\phi_2:= \frac{\varphi - \phi}{(\theta_2 - \theta_0)(\theta^*_0 - \theta^*_2)} -  (\theta_1 - \theta_0)(\theta^*_1 - \theta^*_2).
\end{eqnarray}

\item $E^*_2V$ is an eigenspace for $(A^* - \theta^*_1I)(A - \theta_2I)$ and the corresponding eigenvalue is
\begin{eqnarray}
\label{eq:mu}
\varphi_2:= \frac{\phi - \varphi}{(\theta_0 - \theta_2)(\theta^*_0 - \theta^*_2)} -  (\theta_1 - \theta_2)(\theta^*_1 - \theta^*_2).
\end{eqnarray}
\end{enumerate}

\noindent \rm {\bf Proof}. (i) Let $0 \neq u^*_0 \in E^*_0V$. By row $[0^*D]$ of both Lemma \ref{lem:sixdecomp} and Lemma \ref{lem:actionsplit} the space $E^*_0V$ is an eigenspace for $(A^* - \theta^*_1I)(A - \theta_0I)$; let $\eta$ denote the corresponding eigenvalue. We show $\eta = \varphi_1$. By construction
\begin{eqnarray} \label{eq:eqeta}
(A^* -\theta^*_1I)(A -\theta_0I)u^*_0 = \eta u^*_0.
\end{eqnarray}

\noindent By Lemma \ref{lem:E0stareigen} both
\begin{eqnarray}
\label{eq:eqvarphi}
(A^* -\theta^*_1I)(A^* -\theta^*_2I)(A - \theta_1I)(A - \theta_0I)u^*_0 = \varphi u^*_0,\\
\label{eq:eqphi}
(A^* -\theta^*_1I)(A^* -\theta^*_2I)(A - \theta_1I)(A - \theta_2I)u^*_0 = \phi u^*_0.
\end{eqnarray}
\noindent We subtract (\ref{eq:eqvarphi}) from (\ref{eq:eqphi}) and evaluate the result using (\ref{eq:eqeta}) and $A^*u^*_0 = \theta^*_0u^*_0$; this yields
\begin{eqnarray*}
&&(\phi - \varphi)u^*_0 = (\theta_0 - \theta_2)(A^* -\theta^*_1I)(A^* -\theta^*_2I)(A - \theta_1I)u^*_0\\
&& \qquad \qquad = (\theta_0 - \theta_2)(A^* -\theta^*_1I)(A^* - \theta^*_0I + \theta^*_0I - \theta^*_2I)(A - \theta_0I + \theta_0I - \theta_1I)u^*_0\\
&& \qquad \qquad = (\theta_0 - \theta_2)(\theta^*_0 - \theta^*_2)\biggl(\eta + (\theta^*_0 - \theta^*_1)(\theta_0 - \theta_1)\biggr)u^*_0.
\end{eqnarray*}

\noindent Comparing the coefficients of $u^*_0$ we find $\eta = \varphi_1$ as desired.\\

\noindent (ii) Apply (i) above to $\Phi^\Downarrow$.\\
\noindent (iii) Apply (i) above to $\Phi^\downarrow$ and evaluate the result in light of Lemma \ref{lem:relpararray}(ii).\\
\noindent (iv) Apply (ii) above to $\Phi^\downarrow$ and evaluate the result in light of Lemma \ref{lem:relpararray}(ii). \hfill $\Box$
\end{lemma}

\medskip
\begin{lemma} \label{lem:relphivarphi}
With reference to Notation \ref{note:notephi} and Lemma \ref{lem:releigenseq} we have
\begin{eqnarray*}
\varphi - \varphi_1 \varphi_2 = \phi - \phi_1 \phi_2.
\end{eqnarray*}
\noindent \rm {\bf Proof}. To verify the above equation, eliminate
$\varphi_1$, $\varphi_2$, $\phi_1$ and $\phi_2$ using (\ref{eq:lambda})--(\ref{eq:mu}) and simplify the result. \hfill $\Box$
\end{lemma}

\medskip \begin{lemma} \label{lem:releigenseq2}
With reference to Notation \ref{note:notephi},
\begin{enumerate}
\item $E_0V$ is an eigenspace for $(A - \theta_1I)(A^* - \theta^*_0I)$ and the corresponding eigenvalue is the scalar $\varphi_1$ from (\ref{eq:lambda}).

\item $E_0V$ is an eigenspace for $(A - \theta_1I)(A^* - \theta^*_2I)$ and the corresponding eigenvalue is the scalar $\phi_2$ from (\ref{eq:mu2}).

\item $E_2V$ is an eigenspace for $(A - \theta_1I)(A^* - \theta^*_0I)$ and the corresponding eigenvalue is the scalar $\phi_1$ from (\ref{eq:lambda2}).

\item $E_2V$ is an eigenspace for $(A - \theta_1I)(A^* - \theta^*_2I)$  and the corresponding eigenvalue is the scalar $\varphi_2$ from (\ref{eq:mu}).

\end{enumerate}
\noindent \rm {\bf Proof}. Apply Lemma \ref{lem:releigenseq} to $\Phi^*$ and evaluate the result using (\ref{Phistar}) and Lemma \ref{lem:relpararray}(i). \hfill $\Box$
\end{lemma}

\section{Six bases for $V$}

\noindent In this section we continue to consider the situation of Notation \ref{note:notephi}. Referring to that notation we will define six bases for $V$. The first four will be obtained from the split decompositions of $V$. The other two will consist of an eigenbasis for $A$ and an eigenbasis for $A^*$.\\

\noindent We begin with a definition.

\begin {definition}
\rm \label {def:initvect} With reference to Notation \ref{note:notephi}, we fix a nonzero $\eta^*_0 \in E^*_0V$ and define
\begin{eqnarray*}
\eta_0 &=& (A - \theta_1I)(A - \theta_2I)\eta^*_0, \\
\eta_2 &=& (A - \theta_1I)(A - \theta_0I)\eta^*_0, \\
\eta^*_2 &=& (A^* - \theta^*_1I)(A^* - \theta^*_0I)\eta_2. \\
\end{eqnarray*}
By Lemma \ref{lem:actionsplit} $\eta^*_0 \in E^*_0V$, $\eta^*_2 \in E^*_2V$, $\eta_0 \in E_0V$, $\eta_2 \in E_2V$. By construction and by \cite [Lemma 6.5]{pandq} each of $\eta_0, \eta_2, \eta^*_0, \eta^*_2$ is nonzero.
\end{definition}

\medskip \begin {lemma}
\label {lem:Ubases} With reference to Notation \ref{note:notephi} let $\{U_i\}^2_{i=0}$ denote the decomposition $[0^*2]$. Then the following (i)--(iv) hold.
\begin{enumerate}
\item $\eta^*_0$ is a basis for $U_0.$
\item The vectors $(A - \theta_0I) \eta^*_0,  (A^* - \theta_2^*I) \eta_2$ form a basis for $U_1.$
\item $\eta_2$ is a basis for $U_2.$
\item The sequence
\begin {eqnarray}
\label {eq:Ubases} \eta^*_0, \quad (A - \theta_0I)\eta^*_0, \quad (A^* - \theta_2^*I)\eta_2, \quad \eta_2
\end {eqnarray}
is a basis for $V$.
\end{enumerate}

\noindent \rm {\bf Proof}. By Definition \ref{def:initvect} $\eta^*_0 \in U_0$ and $\eta_2 \in U_2$. Applying row $[0^*D]$ of Lemma \ref{lem:actionsplit}, we find both vectors $(A - \theta_0I) \eta^*_0,  (A^* - \theta_2^*I) \eta_2$ are contained in $U_1$. To finish the proof, we show that these four vectors span $V$. Let $W$ be the subspace of $V$ spanned by these four vectors. We show $W = V$. To do this, we show that $W$ is invariant under the actions of $A$ and $A^*$. Let us examine the actions of $A$ and $A^*$ on these vectors. We begin with $A$. Observe

\begin{eqnarray*}
A \eta^*_0 = \theta_0 \eta^*_0 + (A - \theta_0I)\eta^*_0.
\end{eqnarray*}

\noindent By Definition \ref{def:initvect},
\begin{eqnarray*}
A (A - \theta_0I) \eta^*_0 = \theta_1 (A - \theta_0I) \eta^*_0 + \eta_2.
\end{eqnarray*}

\noindent By Lemma \ref{lem:releigenseq2}(iv),
\begin{eqnarray*}
A (A^* - \theta^*_2I) \eta_2 = \theta_1 (A^* - \theta^*_2I) \eta_2 + \varphi_2 \eta_2.
\end{eqnarray*}

\noindent By (\ref{eq:idem4}),
\begin{eqnarray*}
A \eta_2 = \theta_2 \eta_2.
\end{eqnarray*}

\noindent By these comments $AW \subseteq W$.\\

\noindent Regarding $A^*$ we have the following.
\noindent Since $E^*_0V$ is the eigenspace of $A^*$ corresponding to $\theta^*_0$,
\begin{eqnarray*}
A^* \eta^*_0 = \theta^*_0 \eta^*_0.
\end{eqnarray*}

\noindent By Lemma \ref{lem:releigenseq}(i),
\begin{eqnarray*}
A^* (A - \theta_0I) \eta^*_0 = \varphi_1 \eta^*_0 + \theta^*_1 (A - \theta_0I) \eta^*_0.
\end{eqnarray*}

\noindent By Definition \ref{def:initvect} and Lemma \ref{lem:E0stareigen}(i),
\begin{eqnarray*}
A^* (A^* - \theta^*_2I) \eta_2 = \varphi \eta^*_0 + \theta^*_1(A^* - \theta^*_2I) \eta_2.
\end{eqnarray*}
\noindent Finally,
\begin{eqnarray*}
A^* \eta_2 = (A^* - \theta^*_2I)\eta_2 + \theta^*_2 \eta_2.
\end{eqnarray*}

\noindent By these comments $A^*W \subseteq W$. Observe that $W \neq 0$ since it contains a nonzero vector $\eta^*_0$. Now $W = V$ in view of Definition \ref{def:TDsys}(iv). We have now shown that the vectors (\ref{eq:Ubases}) span $V$ so they form a basis for $V$ and the result follows. \hfill $\Box$
\end{lemma}

\medskip \begin {lemma}
\label {lem:U2bases} With reference to Notation \ref{note:notephi} let $\{U_i\}^2_{i=0}$ denote the decomposition $[0^*0]$.Then the following (i)--(iv) hold.
\begin{enumerate}
\item $\eta^*_0$ is a basis for $U_0.$
\item The vectors $(A - \theta_2I) \eta^*_0,  (A^* - \theta_2^*I) \eta_0$ form a basis for $U_1.$
\item $\eta_0$ is a basis for $U_2.$
\item The sequence
\begin {eqnarray}
\label {eq:U2bases} \eta^*_0, \quad (A - \theta_2I) \eta^*_0, \quad (A^* - \theta_2^*I) \eta_0, \quad \eta_0
\end {eqnarray}
is a basis for $V$.
\end{enumerate}
\noindent \rm {\bf Proof}. Apply Lemma \ref{lem:Ubases} to $\Phi^\Downarrow$ and use Lemma \ref{lem:relpararray}(iii). \hfill $\Box$
\end{lemma}

\medskip \begin {lemma}
\label {lem:U3bases} With reference to Notation \ref{note:notephi} let $\{U_i\}^2_{i=0}$ denote the decomposition $[2^*0]$. Then the following (i)--(iv) hold.
\begin{enumerate}
\item $ \eta^*_2$ is a basis for $U_0.$
\item The vectors $(A - \theta_2I) \eta^*_2,  \varphi(A^* - \theta_0^*I) \eta_0$ form a basis for $U_1$.
\item $\varphi \eta_0$ is a basis for $U_2.$

\item The sequence
\begin {eqnarray}
\label {eq:U3bases} \eta^*_2, \quad (A - \theta_2I) \eta^*_2, \quad \varphi(A^* - \theta^*_0I) \eta_0, \quad \varphi \eta_0
\end {eqnarray}
is a basis for $V$.
\end{enumerate}

\noindent \rm {\bf Proof}. Apply Lemma \ref{lem:U2bases} to $\Phi^\downarrow$ and use Lemma \ref{lem:relpararray}(ii). \hfill $\Box$

\end{lemma}

\medskip \begin {lemma}
\label {lem:U4bases} With reference to Notation \ref{note:notephi} let $\{U_i\}^2_{i=0}$ denote the decomposition $[2^*2]$. Then the following (i)--(iv) hold.
\begin{enumerate}
\item $\eta^*_2$ is a basis for $U_0.$
\item The vectors $(A - \theta_0I) \eta^*_2,  \phi(A^* - \theta_0^*I) \eta_2$ form a basis for $U_1.$
\item $\phi \eta_2$ is a basis for $U_2$.
\item The sequence
\begin {eqnarray}
\label {eq:U4bases}  \eta^*_2, \quad (A - \theta_0I) \eta^*_2, \quad \phi(A^* - \theta_0^*I) \eta_2, \quad \phi \eta_2
\end {eqnarray}
is a basis for $V$.
\end{enumerate}

\noindent \rm {\bf Proof}. Apply Lemma \ref{lem:Ubases} to $\Phi^\downarrow$ and use Lemma \ref{lem:relpararray}(ii). \hfill $\Box$

\end{lemma}

\medskip
\begin{definition} \rm \label{def:defsplitbases}
We refer to the bases (\ref{eq:Ubases})--(\ref{eq:U4bases}) as the {\it split bases} for $V$.
\end{definition}

\noindent We now display an eigenbasis for $A$.

\begin{lemma} \label{lem:Aeigen} With reference to Notation \ref{note:notephi} let $\{U_i\}^2_{i=0}$ denote the decomposition $[02]$. Then the following (i)--(iv) hold.
\begin{enumerate}
\item $\eta_0$ is a basis for $U_0.$
\item The vectors $E_1 \eta^*_0,  E_1 \eta^*_2$ form a basis for $U_1.$
\item $\eta_2$ is a basis for $U_2.$
\item The sequence
\begin {eqnarray}
\label {eq:Abases} \eta_0, \quad E_1 \eta^*_0, \quad E_1 \eta^*_2, \quad \eta_2
\end {eqnarray}
is a basis for $V$.
\end{enumerate}

\noindent \rm {\bf Proof}. (i), (iii) By construction $\eta_0$ is a nonzero vector contained in $U_0$, and $\eta_2$ is a nonzero vector contained in $U_2$. The result follows since $U_0$ and $U_2$ are both $1$-dimensional subspaces of $V$.\\
\noindent (ii) By the comment preceding Definition \ref{def:TDsys}, each of  $E_1 \eta^*_0$, $E_1 \eta^*_2$ is contained in $U_1 = E_1V$. Since $U_1$ has dimension 2, it suffices to show that they are linearly independent. To do this, we write them in terms of the basis (\ref{eq:Ubases}). We claim

\begin {eqnarray}
\label{eq:E1eta0star}
E_1 \eta^*_0 &=& \frac{\eta_2 + (\theta_1 - \theta_2)(A - \theta_0I)\eta^*_0}{(\theta_1 - \theta_0)(\theta_1 - \theta_2)}.
\end{eqnarray}
\noindent To obtain (\ref{eq:E1eta0star}) use the formula for $E_1$ in Notation \ref{note:notephi} and the fact that $A - \theta_2I = A - \theta_1I + (\theta_1 - \theta_2)I$. Then simplify the result using Definition \ref{def:initvect}. We now have (\ref{eq:E1eta0star}). Next we claim

\begin {eqnarray} \label{eq:E1eta2star}
E_1 \eta^*_2 &=&  \frac{\varphi(A - \theta_0I)\eta^*_0}{\theta_1 - \theta_0} + (\theta^*_2 - \theta^*_0)(A^* - \theta^*_2I)\eta_2 + \frac{\varphi +  \varphi_2(\theta_1 - \theta_0)(\theta^*_2 - \theta^*_0)}{(\theta_1 - \theta_0)(\theta_1 - \theta_2)}\eta_2.
\end{eqnarray}

\noindent To obtain (\ref{eq:E1eta2star}) we observe that by row $[0^*D]$ of Lemma \ref{lem:actionsplit},
\begin{eqnarray}
\label{eq:Atheta2eta2}
(A - \theta_2I)\eta_2 = 0.
\end{eqnarray}

\noindent By row $[0D]$ of Lemma \ref{lem:actionsplit},
\begin{eqnarray}
\label{eq:Atheta1theta2theta2star}
(A - \theta_1I)(A - \theta_2I)(A^* - \theta^*_1I)\eta_2 = 0.
\end{eqnarray}

\noindent By Lemma \ref{lem:E2eigen}(i),
\begin{eqnarray}
\label{eq:varphieta2}
(A - \theta_1I)(A - \theta_0I)(A^* - \theta^*_1I)(A^* - \theta^*_2I)\eta_2 = \varphi \eta_2.
\end{eqnarray}

\noindent By Lemma \ref{lem:releigenseq2}(iv),
\begin {eqnarray}
\label{eq:varphi2eta2}
(A - \theta_1I)(A^* - \theta^*_2I)\eta_2 = \varphi_2 \eta_2.
\end{eqnarray}

\noindent By Lemma \ref{lem:E0stareigen}(i) and Definition \ref{def:initvect},
\begin {eqnarray}
\label{eq:Atheta0theta1startheta2star}
(A - \theta_0I)(A^* - \theta^*_1I)(A^* - \theta^*_2I)\eta_2 = \varphi (A - \theta_0I) \eta^*_0.
\end{eqnarray}

\noindent Consider the equation which is $(\theta_1 - \theta_0)^{-1}(\theta_1 - \theta_2)^{-1}$ times (\ref{eq:varphieta2}) plus $(\theta^*_2 - \theta^*_0)(\theta_1 - \theta_2)^{-1}$ times (\ref{eq:varphi2eta2}) plus $(\theta_1 - \theta_0)^{-1}$ times (\ref{eq:Atheta0theta1startheta2star}). Adding $(\theta^*_2 - \theta^*_0)(A^* - \theta^*_2I)\eta_2$ to both sides of this equation and simplifying the result using (\ref{eq:Atheta2eta2}), (\ref{eq:Atheta1theta2theta2star}) and the equation for $E_1$ in Notation \ref{note:notephi}, we routinely obtain (\ref{eq:E1eta2star}).\\\\
\noindent We now compare (\ref{eq:E1eta0star}) and (\ref{eq:E1eta2star}). Observe that the coefficient of $(A^* - \theta^*_2I)\eta_2$ is zero in (\ref{eq:E1eta0star}) and nonzero in (\ref{eq:E1eta2star}). Therefore $E_1 \eta^*_0$ and $E_1 \eta^*_2$ are linearly independent as desired. The result follows.\\
\noindent (iv). Immediate from (i)--(iii) above. \hfill $\Box$
\end{lemma}

\medskip
\noindent Next we display an eigenbasis for $A^*$.

\begin{lemma} \label{lem:Astareigin} With reference to Notation \ref{note:notephi} let $\{U_i\}^2_{i=0}$ denote the decomposition $[0^*2^*]$. Then the following (i)--(iv) hold.
\begin{enumerate}
\item $\eta^*_0$ is a basis for $U_0.$
\item The vectors $E^*_1 \eta_0,  E^*_1 \eta_2$ form a basis for $U_1.$
\item $\eta^*_2$ is a basis for $U_2.$
\item The sequence
\begin {eqnarray}
\label {eq:Astarbases} \eta^*_0, \quad E^*_1 \eta_0, \quad E^*_1 \eta_2, \quad \eta^*_2
\end {eqnarray}
is a basis for $V$.
\end{enumerate}

\medskip
\noindent \rm {\bf Proof}. Apply Lemma \ref{lem:Aeigen} to $\Phi^*$. \hfill $\Box$
\end{lemma}

\medskip
\section{The action of $A$ and $A^*$ on the six bases}

\noindent In this section we display the matrices representing $A$ and $A^*$ with respect to the six bases presented in the previous section. We first recall some basic facts from linear algebra. Let $A$ be a linear transformation on $V$ and let $\lbrace v_i \rbrace_{i=0}^d$ be a basis for $V$. We say that a matrix $B$ {\it represents} $A$ with respect to the basis $\lbrace v_i \rbrace_{i=0}^d$ whenever $Av_j = \sum_{i=0}^d B_{ij} v_i$ for $0 \leq j \leq d$.\\

\noindent We now display the matrices representing $A$ and $A^*$ with respect to the split bases.

\medskip
\begin{theorem}
\label{thm:Asplit}With reference to Notation \ref{note:notephi} the following (i)--(iv) hold.
\begin{enumerate}
\item The matrices representing $A$ and $A^*$ with respect to the basis (\ref {eq:Ubases}) are \\
$$ \left(
\begin{array}{ c c c c }
\theta_0 & 0  &  0    & 0  \\
1 & \theta_1  &  0    & 0  \\
0 & 0  &  \theta_1 &    0  \\
0 & 1  &  \varphi_2    & \theta_2  \\
\end{array}
\right), \qquad
\left(
\begin{array}{ c c c c }
\theta^*_0 & \varphi_1  &  \varphi    & 0  \\
0 & \theta^*_1  &  0   &  0    \\
0  & 0  & \theta^*_1   & 1 \\
0  & 0  & 0  & \theta^*_2 \\
\end{array}
\right)$$ respectively.

\item The matrices representing $A$ and $A^*$ with respect to the basis
(\ref {eq:U2bases}) are \\
$$ \left(
\begin{array}{ c c c c }
\theta_2 & 0  &  0    & 0  \\
1 & \theta_1  &  0    & 0  \\
0 & 0  &  \theta_1   & 0  \\
0 & 1  &  \phi_2    & \theta_0  \\
\end{array}
\right), \qquad \left(
\begin{array}{ c c c c }
\theta^*_0 & \phi_1  &  \phi    & 0  \\
0 & \theta^*_1  &  0   &  0    \\
0  & 0  & \theta^*_1   & 1 \\
0  & 0  & 0  & \theta^*_2 \\
\end{array}
\right)$$ respectively.

\item The matrices representing $A$ and $A^*$ with respect to the basis (\ref {eq:U3bases}) are \\
$$ \left(
\begin{array}{ c c c c }
\theta_2 & 0  &  0    & 0  \\
1 & \theta_1  &  0    & 0  \\
0 & 0  &  \theta_1    & 0  \\
0 & 1  &  \varphi_1    & \theta_0  \\
\end{array}
\right), \qquad \left(
\begin{array}{ c c c c }
\theta^*_2 & \varphi_2  &  \varphi    & 0  \\
0 & \theta^*_1  &  0   &  0    \\
0  & 0  & \theta^*_1   & 1 \\
0  & 0  & 0  & \theta^*_0 \\
\end{array}
\right)$$ respectively.

\item The matrices representing $A$ and $A^*$ with respect to the basis (\ref {eq:U4bases}) are \\
$$ \left(
\begin{array}{ c c c c }
\theta_0 & 0  &  0    & 0  \\
1 & \theta_1  &  0    & 0  \\
0 & 0  &  \theta_1  &  0  \\
0 & 1  &  \phi_1    & \theta_2  \\
\end{array}
\right), \qquad \left(
\begin{array}{ c c c c }
\theta^*_2 & \phi_2  &  \phi    & 0  \\
0 & \theta^*_1  &  0   &  0    \\
0  & 0  & \theta^*_1   & 1 \\
0  & 0  & 0  & \theta^*_0 \\
\end{array}
\right)$$ respectively.
\end{enumerate}

\noindent \rm {\bf Proof}. (i) Immediate from the proof of Lemma \ref{lem:Ubases}.\\
\noindent(ii) Apply (i) above to $\Phi^\Downarrow$ and evaluate the result using Lemmas \ref{lem:E0stareigen}(ii), \ref {lem:relpararray}(iii), \ref{lem:releigenseq}(ii) and \ref{lem:releigenseq2}(iv).\\
\noindent(iii) Apply (i) above to $\Phi^{\downarrow \Downarrow}$ and evaluate the result using Lemmas \ref {lem:relpararray}(ii),(iii), \ref{lem:E2stareigen}(i), \ref{lem:releigenseq}(iv) and \ref{lem:releigenseq2}(i).\\
\noindent(iv) Apply (i) above to $\Phi^\downarrow$ and evaluate the result using Lemmas \ref {lem:relpararray}(ii), \ref{lem:E2stareigen}(ii), \ref{lem:releigenseq}(iii) and \ref{lem:releigenseq2}(iii). \hfill $\Box$
\end{theorem}

\noindent With respect to the eigenbasis of $A$ and $A^*$, we have the following.
\medskip
\begin{theorem}
\label{thm:AAstarbasis}With reference to Notation \ref{note:notephi} the following (i), (ii) hold.
\begin{enumerate}
\item The matrices representing $A$ and $A^*$ with respect to the basis (\ref {eq:Abases}) are \\

\medskip
\centerline{{\rm diag($\theta_0, \theta_1, \theta_1,\theta_2$)},}

$$ {\small \left(
\begin{array}{ c c c c }
\theta^*_0 + \frac{\varphi_1}{\theta_0 - \theta_1} & \frac{\varphi_1}{(\theta_0 - \theta_1)^2(\theta_2 - \theta_0)}  &  \frac{\varphi \phi_2}{(\theta_0 - \theta_1)^2(\theta_2 - \theta_0)}    & 0  \\
\frac{\phi}{\theta^*_0 - \theta^*_2} & \theta^*_1 + \frac{\varphi + \varphi_1(\theta_1 - \theta_2)(\theta^*_0 - \theta^*_2) }{(\theta_1 - \theta_0)(\theta_1 - \theta_2)(\theta^*_0 - \theta^*_2)} &  \frac{\varphi \phi}{(\theta_1 - \theta_0)(\theta_1 - \theta_2)(\theta^*_0 - \theta^*_2)}   &  \frac{\varphi}{\theta^*_0 - \theta^*_2}    \\
\frac{1}{\theta^*_2 - \theta^*_0} &  \frac{1}{(\theta_1 - \theta_0)(\theta_1 - \theta_2)(\theta^*_2 - \theta^*_0)}   & \theta^*_1 + \frac{\varphi + \varphi_2(\theta_1 - \theta_0)(\theta^*_2 - \theta^*_0)}{(\theta_1 - \theta_0)(\theta_1 - \theta_2)(\theta^*_2 - \theta^*_0)}  &  \frac{1}{\theta^*_2 - \theta^*_0}    \\
0  & \frac{\phi_1}{(\theta_1 - \theta_2)^2(\theta_0 - \theta_2)}  & \frac{\phi \varphi_2}{(\theta_1 - \theta_2)^2(\theta_0 - \theta_2)}  & \theta^*_2 + \frac{\varphi_2}{\theta_2 - \theta_1} \\
\end{array}
\right)}$$

respectively.

\item The matrices representing $A$ and $A^*$ with respect to the basis (\ref {eq:Astarbases}) are \\
$$ {\small \left(
\begin{array}{ c c c c }
\theta_0 + \frac{\varphi_1}{\theta^*_0 - \theta^*_1} & \frac{\phi \varphi_1}{(\theta^*_0 - \theta^*_1)^2(\theta^*_2 - \theta^*_0)}  &  \frac{\varphi \phi_1}{(\theta^*_0 - \theta^*_1)^2(\theta^*_2 - \theta^*_0)}    & 0  \\
\frac{1}{\theta_0 - \theta_2} & \theta_1 + \frac{\varphi + \varphi_1(\theta_0 - \theta_2)(\theta^*_1 - \theta^*_2)}{(\theta_0 - \theta_2)(\theta^*_1 - \theta^*_0)(\theta^*_1 - \theta^*_2)} &  \frac{\varphi}{(\theta_0 - \theta_2)(\theta^*_1 - \theta^*_0)(\theta^*_1 - \theta^*_2)}   &  \frac{\varphi}{\theta_0 - \theta_2}    \\
\frac{1}{\theta_2 - \theta_0} &  \frac{\phi}{(\theta_2 - \theta_0)(\theta^*_1 - \theta^*_0)(\theta^*_1 - \theta^*_2)}   & \theta_1 + \frac{\varphi + \varphi_2(\theta_0 - \theta_2)(\theta^*_0 - \theta^*_1) }{(\theta_2 - \theta_0)(\theta^*_1 - \theta^*_0)(\theta^*_1 - \theta^*_2)}  &  \frac{\phi}{\theta_2 - \theta_0}    \\
0  & \frac{\phi_2}{(\theta^*_1 - \theta^*_2)^2(\theta^*_0 - \theta^*_2)}  & \frac{\varphi_2}{(\theta^*_1 - \theta^*_2)^2(\theta^*_0 - \theta^*_2)}  & \theta_2 + \frac{\varphi_2}{\theta^*_2 - \theta^*_1} \\
\end{array}
\right),} $$

\centerline{{\rm diag($\theta^*_0, \theta^*_1, \theta^*_1,\theta^*_2$)},}
respectively.

\end{enumerate}
\noindent {\rm For convenience we will postpone the proof of this theorem until the end of Section $8$.}
\end{theorem}

\medskip
\section{Transition matrices}
\medskip
Suppose we are given two bases for $V$, written
$\lbrace u_i
\rbrace_{i=0}^d$ and
$\lbrace v_i
\rbrace_{i=0}^d$. By the {\it transition matrix}
from
$\lbrace u_i
\rbrace_{i=0}^d$ to
$\lbrace v_i
\rbrace_{i=0}^d$ we mean the matrix
$T$ such that
\begin{eqnarray*}
v_j = \sum_{i=0}^d T_{ij}u_i      \qquad \qquad (0 \leq j\leq d).
\end{eqnarray*}
We recall a few properties of transition matrices.
Let $T$ denote the transition matrix from
$\lbrace u_i
\rbrace_{i=0}^d$ to
$\lbrace v_i
\rbrace_{i=0}^d$. Then $T^{-1}$ exists and equals
the transition matrix from
$\lbrace v_i
\rbrace_{i=0}^d$ to
$\lbrace u_i
\rbrace_{i=0}^d$.
Let $\lbrace w_i \rbrace_{i=0}^d$ denote a basis for $V$, and let
$S$ denote the transition matrix from
$\lbrace v_i
\rbrace_{i=0}^d$ to
$\lbrace w_i \rbrace_{i=0}^d$. Then $TS$ is the transition
matrix from $\lbrace u_i \rbrace_{i=0}^d$  to
$\lbrace w_i \rbrace_{i=0}^d$.\\

\noindent We also recall how the transition matrices and the matrices representing a linear transformation are related. Let $A$ be a linear transformation on $V$ and let $B$ denote the matrix that represents $A$ with respect to $\lbrace u_i
\rbrace_{i=0}^d$. Then the matrix that represents $A$ with respect to $\lbrace v_i
\rbrace_{i=0}^d$ is given by $T^{-1}BT$.\\

\medskip
\noindent For every ordered pair of bases among (\ref{eq:Ubases})--(\ref{eq:Abases}) and (\ref{eq:Astarbases}), we now examine the transition matrices.\\

\medskip
\noindent In the next two theorems we display the transition matrices from one split basis to another.\\

\begin{theorem} \label {thm:TransMat}
With reference to Notation \ref{note:notephi} the following (i)--(iv) hold.
\begin{enumerate}
\item The transition matrix from the basis (\ref{eq:Ubases})  to the basis (\ref{eq:U2bases}) is \\
$$ \left(
\begin{array}{ c c c c }
1 & \theta_0 - \theta_2 & (\theta_0 - \theta_2) \phi_2  &  (\theta_0 - \theta_2)(\theta_0 - \theta_1)  \\
0 & 1  &  (\theta_0 - \theta_2)(\theta^*_1 - \theta^*_2)    & \theta_0 - \theta_2  \\
0 & 0  &  1   & 0  \\
0 & 0  &  0    & 1  \\
\end{array}
\right),
$$

and the transition matrix from the basis  (\ref{eq:U2bases}) to the basis (\ref{eq:Ubases}) is \\
$$ \left(
\begin{array}{ c c c c }
1 & \theta_2 - \theta_0 & (\theta_2 - \theta_0) \varphi_2  &  (\theta_2 - \theta_0)(\theta_2 - \theta_1)  \\
0 & 1  &  (\theta_2 - \theta_0)(\theta^*_1 - \theta^*_2)    & \theta_2 - \theta_0  \\
0 & 0  &  1   & 0  \\
0 & 0  &  0    & 1  \\
\end{array}
\right).
$$

\item The transition matrix from the basis  (\ref{eq:U2bases}) to the basis (\ref{eq:U3bases}) is \\
$$ \left(
\begin{array}{ c c c c }
\phi  & 0  &  0   & 0  \\
0 & \phi  &  0    & 0  \\
\theta^*_2 - \theta^*_0  & (\theta^*_2 - \theta^*_0)(\theta_1 - \theta_2)  &  \varphi   & 0  \\
(\theta^*_2 - \theta^*_0)(\theta^*_2 - \theta^*_1) &(\theta^*_2 - \theta^*_0)\varphi_2 & (\theta^*_2 - \theta^*_0)\varphi  &  \varphi \\
\end{array}
\right),
$$

and the transition matrix from the basis  (\ref{eq:U3bases}) to the basis (\ref{eq:U2bases}) is \\
$$ \left(
\begin{array}{ c c c c }
\phi^{-1} & 0 & 0  & 0  \\
0 & \phi^{-1}  &  0    & 0  \\
(\theta^*_0 - \theta^*_2)\varphi^{-1} \phi^{-1}  & (\theta^*_0 - \theta^*_2)(\theta_1 - \theta_2)\varphi^{-1} \phi^{-1}  &  \varphi^{-1}   & 0  \\
(\theta^*_0 - \theta^*_2)(\theta^*_0 - \theta^*_1)\varphi^{-1} \phi^{-1} & (\theta^*_0 - \theta^*_2)\phi_1\varphi^{-1} \phi^{-1}  &  (\theta^*_0 - \theta^*_2)\varphi^{-1}    &  \varphi^{-1} \\
\end{array}
\right).
$$

\item The transition matrix from the basis (\ref{eq:U3bases}) to the basis (\ref{eq:U4bases}) is \\
$$ \left(
\begin{array}{ c c c c }
1 & \theta_2 - \theta_0 & (\theta_2 - \theta_0) \phi_1  &  (\theta_2 - \theta_0)(\theta_2 - \theta_1)  \\
0 & 1  &  (\theta_2 - \theta_0)(\theta^*_1 - \theta^*_0)    & \theta_2 - \theta_0  \\
0 & 0  &  1   & 0  \\
0 & 0  &  0    & 1  \\
\end{array}
\right),
$$

and the transition matrix from the basis (\ref{eq:U4bases}) to the basis (\ref{eq:U3bases}) is \\
$$ \left(
\begin{array}{ c c c c }
1 & \theta_0 - \theta_2 & (\theta_0 - \theta_2) \varphi_1  &  (\theta_0 - \theta_2)(\theta_0 - \theta_1)  \\
0 & 1  &  (\theta_0 - \theta_2)(\theta^*_1 - \theta^*_0)    & \theta_0 - \theta_2  \\
0 & 0  &  1   & 0  \\
0 & 0  &  0    & 1  \\
\end{array}
\right).
$$

\item The transition matrix from the basis (\ref{eq:U4bases}) to the basis (\ref{eq:Ubases}) is \\
$$ \left(
\begin{array}{ c c c c }
\varphi^{-1} & 0 & 0  & 0  \\
0 & \varphi^{-1}  &  0    & 0  \\
(\theta^*_0 - \theta^*_2)\varphi^{-1} \phi^{-1}  & (\theta^*_0 - \theta^*_2)(\theta_1 - \theta_0)\varphi^{-1} \phi^{-1}  &  \phi^{-1}   & 0  \\
(\theta^*_0 - \theta^*_2)(\theta^*_0 - \theta^*_1)\varphi^{-1} \phi^{-1} & (\theta^*_0 - \theta^*_2)\varphi_1\varphi^{-1} \phi^{-1}  &  (\theta^*_0 - \theta^*_2) \phi^{-1}    & \phi^{-1}  \\
\end{array}
\right),
$$

and the transition matrix from the basis (\ref{eq:Ubases}) to the basis (\ref{eq:U4bases}) is \\
$$ \left(
\begin{array}{ c c c c }
\varphi & 0 & 0  & 0  \\
0 & \varphi  &  0    & 0  \\
\theta^*_2 - \theta^*_0  & (\theta^*_2 - \theta^*_0)(\theta_1 - \theta_0)  &  \phi   & 0  \\
(\theta^*_2 - \theta^*_0)(\theta^*_2 - \theta^*_1) & (\theta^*_2 - \theta^*_0)\phi_2  &  (\theta^*_2 - \theta^*_0)\phi & \phi \\
\end{array}
\right).
$$
\end{enumerate}

\medskip
\noindent \rm {\bf Proof}. (i) We first obtain the transition matrix from the basis (\ref{eq:Ubases}) to the basis (\ref{eq:U2bases}). To get the first column of this matrix, note that the first basis vectors in (\ref{eq:Ubases}) and (\ref{eq:U2bases}) are the same. To obtain the second column, observe
\begin{eqnarray*}
(A - \theta_2I)\eta^*_0 &=& (\theta_0 - \theta_2)\eta^*_0 + (A - \theta_0I)\eta^*_0.
\end{eqnarray*}

\noindent To obtain the third column we observe the following.
\noindent By Definition \ref{def:initvect},
\begin{eqnarray}
\label{eq:eta0toeta2}
\eta_0 = \eta_2 + (\theta_0 - \theta_2)(A - \theta_0I)\eta^*_0 + (\theta_0 - \theta_2)(\theta_0 - \theta_1)\eta^*_0.
\end{eqnarray}

\noindent By Lemma \ref{lem:releigenseq}(i),
\begin{eqnarray}
\label{eq:AstarA0eta0star}
(A^* - \theta^*_2I)(A - \theta_0I)\eta^*_0 = \varphi_1 \eta^*_0 + (\theta^*_1 - \theta^*_2)(A - \theta_0I)\eta^*_0.
\end{eqnarray}

\noindent By row $[0^*D]$ of Lemma \ref{lem:actionsplit} $(A^* - \theta^*_0I)\eta^*_0 = 0$ so
\begin{eqnarray}
\label{Astartheta2stareta0star}
(A^* - \theta^*_2I)\eta^*_0 = (\theta^*_0 - \theta^*_2)\eta^*_0.
\end{eqnarray}

\noindent By (\ref{eq:lambda}), (\ref{eq:mu2}),
\begin{eqnarray}
\label{eq:varphi1tophi2}
\varphi_1 + (\theta_0 - \theta_1)(\theta^*_0 - \theta^*_2) = \phi_2.
\end{eqnarray}

\noindent Applying $A^* - \theta^*_2I$ to both sides of (\ref{eq:eta0toeta2}) and simplifying the result using (\ref{eq:AstarA0eta0star})--(\ref{eq:varphi1tophi2}), we obtain
\begin{eqnarray*}
(A^* - \theta^*_2I)\eta_0 = (\theta_0 - \theta_2)\phi_2 \eta^*_0 +
(\theta_0 - \theta_2)(\theta^*_1 - \theta^*_2)(A -
\theta_0I)\eta^*_0 + (A^* - \theta^*_2I)\eta_2
\end{eqnarray*}
\noindent as desired.\\

\medskip
\noindent The fourth column follows from (\ref{eq:eta0toeta2}).\\
\medskip

\noindent We have now obtained the transition matrix from the basis (\ref{eq:Ubases}) to the basis (\ref{eq:U2bases}). The transition matrix from the basis (\ref{eq:U2bases}) to the basis (\ref{eq:Ubases}) is as shown, since it is routine to verify that the product of this matrix and the previous matrix is the identity.\\

\noindent (ii) We first obtain the transition matrix from the basis (\ref{eq:U2bases}) to the basis (\ref{eq:U3bases}). To obtain the first column, pick $0 \neq \eta^*_0 \in E_0V$. By Definition \ref{def:initvect}, both
\begin{eqnarray}
\label{eq:eta2startoeta2}
\eta^*_2 &=& (A - \theta^*_2I)(A - \theta^*_0I)\eta_2,\\
\label{eq:eta2toeta0eta0star}
\eta_2 &=& \eta_0 + (\theta_2 - \theta_0)(A - \theta_1I)\eta^*_0.
\end{eqnarray}

\noindent By Lemma \ref{lem:E0stareigen} and Definition \ref{def:initvect},
\begin{eqnarray}
\label{eq:phieta0star}
(A^* - \theta^*_1I)(A^* - \theta^*_2I)\eta_0 = \phi\eta^*_0.
\end{eqnarray}

\noindent By row $[0^*D^*]$ of Lemma \ref{lem:actionsplit},
\begin{eqnarray}
\label{eq:zerostarDstareta0star}
(A^* - \theta^*_1I)(A^* - \theta^*_0I)(A - \theta_1I)\eta^*_0 = 0.
\end{eqnarray}

\noindent Eliminate $\eta_2$ in (\ref{eq:eta2startoeta2}) using (\ref{eq:eta2toeta0eta0star}) and simplify the result using (\ref{eq:phieta0star}), (\ref{eq:zerostarDstareta0star}) to obtain
\begin{eqnarray}
\label{eq:firstcol22to23}
\eta^*_2 = \phi \eta^*_0 + (\theta^*_2 - \theta^*_0)(A^* - \theta^*_2I)\eta_0 + (\theta^*_2 - \theta^*_0)(\theta^*_2 - \theta^*_1)\eta_0
\end{eqnarray}

\noindent as desired.\\

\medskip
\noindent To obtain the second column, observe that by row $[0^*0]$ of Lemma \ref{lem:actionsplit},
\begin{eqnarray}
\label{eq:0star0eta0}
(A - \theta_0I)\eta_0 = 0.
\end{eqnarray}

\noindent By Lemma \ref{lem:releigenseq2}(ii),
\begin{eqnarray}
\label{eq:varphi2eta0}
(A - \theta_1I)(A^* - \theta^*_2I)\eta_0 = \phi_2 \eta_0.
\end{eqnarray}

\noindent Note that both
\begin{eqnarray}
\label{eq:derive1}
(A - \theta_2I)(A^* - \theta^*_2I)\eta_0 &=& (\theta_1 - \theta_2)(A^* - \theta^*_2I)\eta_0 + (A - \theta_1I)(A^* - \theta^*_2I)\eta_0,\\
\label{eq:derive2}
(A - \theta_2I)\eta_0 &=& (\theta_0 - \theta_2)\eta_0 + (A - \theta_0I)\eta_0.
\end{eqnarray}

\noindent By (\ref{eq:mu2}), (\ref{eq:mu}),
\begin{eqnarray}
\label{eq:phi2tovarphi2}
\phi_2 + (\theta^*_2 - \theta^*_1)(\theta_0 - \theta_2) = \varphi_2.
\end{eqnarray}

\noindent Applying $A - \theta_2I$ to (\ref{eq:firstcol22to23}) and simplifying the result using (\ref{eq:0star0eta0})--(\ref{eq:phi2tovarphi2}), we obtain the second column.\\
\medskip

\noindent The third column follows from the fact that
\begin{eqnarray*}
(A^* - \theta^*_0I)\varphi \eta_0 = \varphi(A^* - \theta^*_2I)\eta_0 + (\theta^*_2 - \theta^*_0)\varphi \eta_0.
\end{eqnarray*}

\noindent To obtain the fourth column, note that the fourth basis vectors in (\ref{eq:U2bases}) and (\ref{eq:U3bases}) are the same. We have now obtained the transition matrix from the basis (\ref{eq:U2bases}) to the basis (\ref{eq:U3bases}). \\\\
The transition matrix from the basis (\ref{eq:U3bases}) to the basis (\ref{eq:U2bases}) is as shown, since it is routine to verify that the product of this matrix and the previous matrix is the identity.\\

\noindent (iii) Apply (i) above to $\Phi^{\downarrow \Downarrow}$ and use Lemmas \ref{lem:relpararray}(ii),(iii) and \ref{lem:releigenseq}.\\
\noindent (iv) Apply (ii) above to $\Phi^\Downarrow$ and use Lemmas \ref{lem:relpararray}(iii) and \ref{lem:releigenseq}. \hfill $\Box$
\end{theorem}

\medskip
\begin{theorem} \label {thm:TransMat2}
With reference to Notation \ref{note:notephi} the following (i), (ii) hold.
\begin{enumerate}
\item The transition matrix from the basis (\ref{eq:Ubases}) to the basis (\ref{eq:U3bases}) is \\
{\footnotesize $$ \left(
\begin{array}{ c c c c }
\varphi &  (\theta_0 - \theta_2)\varphi & (\theta_0 - \theta_2)\varphi \varphi_1 & (\theta_0 - \theta_2)(\theta_0 - \theta_1)\varphi\\
0 & \varphi & (\theta_0 - \theta_2)(\theta^*_1 - \theta^*_0)\varphi & (\theta_0 - \theta_2)\varphi  \\
\theta^*_2 - \theta^*_0  & (\theta^*_2 - \theta^*_0)(\theta_1 - \theta_2)  &  \varphi   & 0  \\
(\theta^*_2 - \theta^*_0)(\theta^*_2 - \theta^*_1) & (\theta^*_2 - \theta^*_0)\varphi_2  &  (\theta^*_2 - \theta^*_0)\varphi    &  \varphi \\
\end{array}
\right),
$$}

and the transition matrix from the basis  (\ref{eq:U3bases}) to the basis (\ref{eq:Ubases}) is \\
{\footnotesize $$ \left(
\begin{array}{ c c c c }
\phi^{-1}  & (\theta_2 - \theta_0)\phi^{-1}  &  (\theta_2 - \theta_0) \varphi_2 \phi^{-1}   & (\theta_2 - \theta_0)(\theta_2 - \theta_1)\phi^{-1}  \\
0 & \phi^{-1}  &  (\theta_2 - \theta_0)(\theta^*_1 - \theta^*_2)\phi^{-1}    & (\theta_2 - \theta_0)\phi^{-1}  \\
(\theta^*_0 - \theta^*_2)\varphi^{-1} \phi^{-1}  & (\theta^*_0 - \theta^*_2)(\theta_1 - \theta_0)\varphi^{-1} \phi^{-1}  &  \phi^{-1}   & 0  \\
(\theta^*_0 - \theta^*_2)(\theta^*_0 - \theta^*_1)\varphi^{-1} \phi^{-1} & (\theta^*_0 -\theta^*_2)\varphi_1\varphi^{-1} \phi^{-1} & (\theta^*_0 - \theta^*_2)\phi^{-1} & \phi^{-1}\\
\end{array}
\right).
$$}

\item The transition matrix from the basis (\ref{eq:U2bases}) to the basis (\ref{eq:U4bases}) is \\
{\footnotesize $$  \left(
\begin{array}{ c c c c }
\phi  & (\theta_2 - \theta_0)\phi  &  (\theta_2 - \theta_0)\phi \phi_1   & (\theta_2 - \theta_0)(\theta_2 - \theta_1)\phi  \\
0 & \phi  &  (\theta_2 - \theta_0)(\theta^*_1 - \theta^*_0)\phi    & (\theta_2 - \theta_0)\phi  \\
\theta^*_2 - \theta^*_0  & (\theta^*_2 - \theta^*_0)(\theta_1 - \theta_0)  &  \phi   & 0  \\
(\theta^*_2 - \theta^*_0)(\theta^*_2 - \theta^*_1) & (\theta^*_2 -\theta^*_0)\phi_2 & (\theta^*_2 - \theta^*_0)\phi & \phi\\
\end{array}
\right),
$$}

and the transition matrix from the basis (\ref{eq:U4bases}) to the basis (\ref{eq:U2bases}) is \\
{\footnotesize $$  \left(
\begin{array}{ c c c c }
\varphi^{-1} & (\theta_0 - \theta_2)\varphi^{-1} & (\theta_0 - \theta_2)\phi_2\varphi^{-1} & (\theta_0 - \theta_2)(\theta_0 - \theta_1)\varphi^{-1}\\
0 & \varphi^{-1} & (\theta_0 - \theta_2)(\theta^*_1 - \theta^*_2)\varphi^{-1} & (\theta_0 - \theta_2)\varphi^{-1}  \\
(\theta^*_0 - \theta^*_2)\varphi^{-1} \phi^{-1}  & (\theta^*_0 - \theta^*_2)(\theta_1 - \theta_2)\varphi^{-1} \phi^{-1}  &  \varphi^{-1}   & 0  \\
(\theta^*_0 - \theta^*_2)(\theta^*_0 - \theta^*_1)\varphi^{-1} \phi^{-1} & (\theta^*_0 - \theta^*_2) \phi_1\varphi^{-1} \phi^{-1}  &  (\theta^*_0 - \theta^*_2)\varphi^{-1}    &  \varphi^{-1} \\
\end{array}
\right).
$$}
\end{enumerate}

\medskip
\noindent \rm {\bf Proof}. (i) To obtain the transition matrix from the basis (\ref{eq:Ubases}) to the basis (\ref{eq:U3bases}), compute the product of the transition matrix from the basis (\ref{eq:Ubases}) to the basis (\ref{eq:U2bases}) given in Theorem \ref{thm:TransMat}(i) and the transition matrix from the basis (\ref{eq:U2bases}) to the basis (\ref{eq:U3bases}) given in Theorem \ref{thm:TransMat}(ii). Simplify the product using (\ref{eq:lambda})--(\ref{eq:mu}). To obtain the transition matrix from the basis (\ref{eq:U3bases}) to the basis (\ref{eq:Ubases}), compute the product of the transition matrix from the basis (\ref{eq:U3bases}) to the basis (\ref{eq:U2bases}) given in Theorem \ref{thm:TransMat}(ii) and the transition matrix from the basis (\ref{eq:U2bases}) to the basis (\ref{eq:Ubases}) given in Theorem \ref{thm:TransMat}(i). Simplify the product using (\ref{eq:lambda})--(\ref{eq:mu}).\\\\ 
\noindent (ii) To obtain the transition matrix from the basis (\ref{eq:U2bases}) to the basis (\ref{eq:U4bases}), compute the product of the transition matrix from the basis (\ref{eq:U2bases}) to the basis (\ref{eq:U3bases}) given in Theorem \ref{thm:TransMat}(ii) and the transition matrix from the basis (\ref{eq:U3bases}) to the basis (\ref{eq:U4bases}) given in Theorem \ref{thm:TransMat}(iii). Simplify the product using (\ref{eq:lambda})--(\ref{eq:mu}). To obtain the transition matrix from the basis (\ref{eq:U4bases}) to the basis (\ref{eq:U2bases}), compute the product of the transition matrix from the basis (\ref{eq:U4bases}) to the basis (\ref{eq:U3bases}) given in Theorem \ref{thm:TransMat}(iii) and the transition matrix from the basis (\ref{eq:U3bases}) to the basis (\ref{eq:U2bases}) given in Theorem \ref{thm:TransMat}(ii). Simplify the product using (\ref{eq:lambda})--(\ref{eq:mu}). \hfill $\Box$
\end{theorem}

\medskip
\noindent In the next theorem we display the the transition matrix between a split basis and  an eigenbasis for $A$.\\
\medskip
\begin{theorem} \label {thm:TransMat3}
With reference to Notation \ref{note:notephi} the following (i)--(iv) hold.
\begin{enumerate}
\item The transition matrix from the basis (\ref{eq:Ubases})  to the basis (\ref{eq:Abases}) is \\
$$ \left(
\begin{array}{ c c c c }
(\theta_0 - \theta_1)(\theta_0 - \theta_2) & 0 & 0  &  0  \\
\theta_0 - \theta_2 & (\theta_1 - \theta_0)^{-1}  &  \varphi(\theta_1 - \theta_0)^{-1}& 0  \\
0 & 0  &  \theta^*_2 - \theta^*_0   & 0  \\
1 & \frac{1}{(\theta_1 - \theta_0)(\theta_1-\theta_2)}  &  \frac{\varphi + \varphi_2(\theta_1 - \theta_0)(\theta^*_2-\theta^*_0)}{(\theta_1 - \theta_0)(\theta_1-\theta_2)}    & 1  \\
\end{array}
\right),
$$

and the transition matrix from the basis  (\ref{eq:Abases}) to the basis (\ref{eq:Ubases}) is \\
$$ \left(
\begin{array}{ c c c c }
\frac{1}{(\theta_0 - \theta_1)(\theta_0 - \theta_2)} & 0 & 0  &  0  \\
1 & \theta_1 - \theta_0 &  \varphi(\theta^*_0 - \theta^*_2)^{-1} & 0  \\
0 & 0  &  (\theta^*_2 - \theta^*_0)^{-1}   & 0  \\
\frac{1}{(\theta_2 - \theta_0)(\theta_2 - \theta_1)} & \frac{1}{\theta_2 - \theta_1}  &  \frac{\varphi_2}{\theta_2-\theta_1}    & 1  \\
\end{array}
\right).
$$

\item The transition matrix from the basis  (\ref{eq:U2bases}) to the basis (\ref{eq:Abases}) is \\
$$ \left(
\begin{array}{ c c c c }
0  & 0  &  0   & (\theta_2 - \theta_0)(\theta_2 - \theta_1)  \\
0 & (\theta_1 - \theta_2)^{-1}  &  \phi (\theta_1 - \theta_2)^{-1}    & \theta_2 - \theta_0  \\
0  & 0  &  \theta^*_2 - \theta^*_0   & 0  \\
1 & \frac{1}{(\theta_1 - \theta_0)(\theta_1 - \theta_2)} & \frac{\phi + \phi_2(\theta_1 - \theta_2)(\theta^*_2 - \theta^*_0)}{(\theta_1 - \theta_0)(\theta_1 - \theta_2)}      & 1 \\
\end{array}
\right),
$$

and the transition matrix from the basis  (\ref{eq:Abases}) to the basis (\ref{eq:U2bases}) is \\
$$ \left(
\begin{array}{ c c c c }
\frac{1}{(\theta_0 - \theta_1)(\theta_0 - \theta_2)} & \frac{1}{\theta_0 - \theta_1} & \frac{\phi_2}{\theta_0 - \theta_1}  &  1  \\
1 & \theta_1 - \theta_2 &   \phi(\theta^*_0 - \theta^*_2)^{-1} & 0  \\
0 & 0  &  (\theta^*_2 - \theta^*_0)^{-1}   & 0  \\
\frac{1}{(\theta_2 - \theta_0)(\theta_2 - \theta_1)} & 0  &  0    & 0  \\
\end{array}
\right).
$$

\item The transition matrix from the basis (\ref{eq:U3bases}) to the basis (\ref{eq:Abases}) is \\
$$ \left(
\begin{array}{ c c c c }
0  & 0  &  0   & (\theta_2 - \theta_0)(\theta_2 - \theta_1)\phi^{-1}  \\
0 & (\theta_1 - \theta_2)^{-1}\phi^{-1}  &  (\theta_1 - \theta_2)^{-1}   &  (\theta_2 - \theta_0)\phi^{-1}  \\
0  & (\theta^*_0 - \theta^*_2)\varphi^{-1}\phi^{-1}  &  0   & 0  \\
\varphi^{-1}  & \frac{\phi + \phi_1(\theta_1 - \theta_0)(\theta^*_0 - \theta^*_2)}{(\theta_1 - \theta_0)(\theta_1 - \theta_2)\varphi \phi} & \frac{1}{(\theta_1 - \theta_0)(\theta_1 - \theta_2)}      & \phi^{-1} \\
\end{array}
\right),
$$

and the transition matrix from the basis (\ref{eq:Abases}) to the basis (\ref{eq:U3bases}) is \\
$$ \left(
\begin{array}{ c c c c }
\frac{\varphi}{(\theta_0 - \theta_1)(\theta_0 - \theta_2)} & \frac{\varphi}{\theta_0 - \theta_1} & \frac{\varphi\varphi_1}{\theta_0 - \theta_1}  &  \varphi  \\
0 & 0  &  \varphi\phi(\theta^*_0 - \theta^*_2)^{-1}    & 0  \\
1 & \theta_1 - \theta_2  &  \varphi(\theta^*_2 - \theta^*_0)^{-1}   & 0  \\
\frac{\phi}{(\theta_0 - \theta_2)(\theta_1 - \theta_2)} & 0  &  0    & 0  \\
\end{array}
\right).
$$

\item The transition matrix from the basis (\ref{eq:U4bases}) to the basis (\ref{eq:Abases}) is \\
$$ \left(
\begin{array}{ c c c c }
(\theta_0 - \theta_1)(\theta_0 - \theta_2)\varphi^{-1} & 0 & 0  &  0  \\
(\theta_0 - \theta_2)\varphi^{-1} & \varphi^{-1}(\theta_1 - \theta_0)^{-1}  &  (\theta_1 - \theta_0)^{-1} & 0  \\
0  &  (\theta^*_0 - \theta^*_2)\varphi^{-1}\phi^{-1} & 0   & 0  \\
\varphi^{-1} &  \frac{\varphi + \varphi_1(\theta_1 - \theta_2)(\theta^*_0-\theta^*_2)}{(\theta_1 - \theta_0)(\theta_1-\theta_2)\varphi \phi}  & \frac{1}{(\theta_1 - \theta_0)(\theta_1-\theta_2)} & \phi^{-1}\\
\end{array}
\right),
$$
and the transition matrix from the basis (\ref{eq:Abases}) to the basis (\ref{eq:U4bases}) is \\
$$ \left(
\begin{array}{ c c c c }
\frac{\varphi}{(\theta_0 - \theta_1)(\theta_0 - \theta_2)} & 0 & 0  &  0  \\
0 & 0  &  \varphi\phi(\theta^*_0 - \theta^*_2)^{-1}   & 0  \\
1 & \theta_1 - \theta_0  &  \phi(\theta^*_2 - \theta^*_0)^{-1} & 0  \\
\frac{\phi}{(\theta_0 - \theta_2)(\theta_1 - \theta_2)} & \frac{\phi}{\theta_2 - \theta_1}  &  \frac{\phi\phi_1}{\theta_2-\theta_1}    & \phi  \\
\end{array}
\right).
$$
\end{enumerate}

\medskip
\noindent \rm {\bf Proof}. (i) We first obtain the transition matrix from the basis (\ref{eq:Ubases}) to the basis (\ref{eq:Abases}). To obtain the first column, we observe that by Definition \ref{def:initvect},
\begin{eqnarray*}
\eta_0 &=& (A - \theta_1I)(A - \theta_2I)\eta^*_0\\
&=& (A - \theta_1I)(A - \theta_0I)\eta^*_0 + (\theta_0 - \theta_2)(A - \theta_1I)\eta^*_0\\
&=& \eta_2 + (\theta_0 - \theta_2)(A - \theta_0I)\eta^*_0 + (\theta_0 - \theta_1)(\theta_0 - \theta_2)\eta^*_0.
\end{eqnarray*}

\noindent The second and third columns follow from (\ref{eq:E1eta0star}) and (\ref{eq:E1eta2star}) respectively.\\

\noindent The last column follows since the fourth basis vectors in (\ref{eq:Ubases}) and (\ref{eq:Abases}) are the same. We have now obtained the transition matrix from (\ref{eq:Ubases}) to (\ref{eq:Abases}).\\

\medskip
\noindent The transition matrix from the basis  (\ref{eq:Abases}) to the basis (\ref{eq:Ubases}) is as shown, since it is routine to verify that the  product of this matrix and the previous matrix is the identity.\\\\
\noindent (ii) To obtain the transition matrix from the basis (\ref{eq:U2bases}) to the basis (\ref{eq:Abases}), compute the product of the transition matrix from the basis (\ref{eq:U2bases}) to the basis (\ref{eq:Ubases}) given in Theorem \ref{thm:TransMat}(i) and the transition matrix from the basis (\ref{eq:Ubases}) to the basis (\ref{eq:Abases}) given in Theorem \ref{thm:TransMat3}(i). Simplify the product using (\ref{eq:lambda})--(\ref{eq:mu}). To obtain the transition matrix from the basis (\ref{eq:Abases}) to the basis (\ref{eq:U2bases}), compute the product of the transition matrix from the basis (\ref{eq:Abases}) to the basis (\ref{eq:Ubases}) given in Theorem \ref{thm:TransMat3}(i) and the transition matrix from the basis (\ref{eq:Ubases}) to the basis (\ref{eq:U2bases}) given in Theorem \ref{thm:TransMat}(i). Simplify the product using (\ref{eq:lambda})--(\ref{eq:mu}).\\\\
\noindent (iii) To obtain the transition matrix from the basis (\ref{eq:U3bases}) to the basis (\ref{eq:Abases}), compute the product of the transition matrix from the basis (\ref{eq:U3bases}) to the basis (\ref{eq:Ubases}) given in Theorem \ref{thm:TransMat2}(i) and the transition matrix from the basis (\ref{eq:Ubases}) to the basis (\ref{eq:Abases}) given in Theorem \ref{thm:TransMat3}(i). Simplify the product using (\ref{eq:lambda})--(\ref{eq:mu}). To obtain the transition matrix from the basis (\ref{eq:Abases}) to the basis (\ref{eq:U3bases}), compute the product of the transition matrix from the basis (\ref{eq:Abases}) to the basis (\ref{eq:Ubases}) given in Theorem \ref{thm:TransMat3}(i) and the transition matrix from the basis (\ref{eq:Ubases}) to the basis (\ref{eq:U3bases}) given in Theorem \ref{thm:TransMat2}(i). Simplify the product using (\ref{eq:lambda})--(\ref{eq:mu}).\\\\
\noindent (iv) To obtain the transition matrix from the basis (\ref{eq:U4bases}) to the basis (\ref{eq:Abases}), compute the product of the transition matrix from the basis (\ref{eq:U4bases}) to the basis (\ref{eq:Ubases}) given in Theorem \ref{thm:TransMat}(iv) and the transition matrix from the basis (\ref{eq:Ubases}) to the basis (\ref{eq:Abases}) given in Theorem \ref{thm:TransMat3}(i). Simplify the product using (\ref{eq:lambda})--(\ref{eq:mu}). To obtain the transition matrix from the basis (\ref{eq:Abases}) to the basis (\ref{eq:U4bases}), compute the product of the transition matrix from the basis (\ref{eq:Abases}) to the basis (\ref{eq:Ubases}) given in Theorem \ref{thm:TransMat3}(i) and the transition matrix from the basis (\ref{eq:Ubases}) to the basis (\ref{eq:U4bases}) given in Theorem \ref{thm:TransMat}(iv). Simplify the product using (\ref{eq:lambda})--(\ref{eq:mu}). \hfill $\Box$
\end{theorem}

\bigskip
\noindent In the next theorem we display the the transition matrix between a split basis and  an eigenbasis for $A^*$.\\
\medskip
\begin{theorem} \label {thm:TransMat4}
With reference to Notation \ref{note:notephi} the following (i)--(iv) hold.
\begin{enumerate}
\item The transition matrix from the basis (\ref{eq:Ubases}) to the basis (\ref{eq:Astarbases}) is \\
$$ \left(
\begin{array}{ c c c c }
1 &  \frac{\phi + \phi_2(\theta_0 - \theta_2)(\theta^*_1 - \theta^*_0)}{(\theta^*_1 - \theta^*_0)(\theta^*_1 - \theta^*_2)} & \frac{\varphi}{(\theta^*_1 - \theta^*_0)(\theta^*_1 - \theta^*_2)} & \varphi\\
0 & \theta_0 - \theta_2 & 0 & 0  \\
0  & (\theta^*_1 - \theta^*_2)^{-1}  &  (\theta^*_1 - \theta^*_2)^{-1}   & \theta^*_2 - \theta^*_0  \\
0 & 0  &  0    &  (\theta^*_2 - \theta^*_0)(\theta^*_2 - \theta^*_1) \\
\end{array}
\right),
$$

and the transition matrix from the basis  (\ref{eq:Astarbases}) to the basis (\ref{eq:Ubases}) is \\
$$ \left(
\begin{array}{ c c c c }
1 &  \frac{\varphi_1}{\theta^*_0 - \theta^*_1} & \frac{\varphi}{\theta^*_0 - \theta^*_1} & \frac{\varphi}{(\theta^*_0 - \theta^*_1)(\theta^*_0 - \theta^*_2)}\\
0 & (\theta_0 - \theta_2)^{-1} & 0 & 0  \\
0  & (\theta_2 - \theta_0)^{-1}  &  \theta^*_1 - \theta^*_2   &  1 \\
0 & 0  &  0    & \frac{1}{(\theta^*_0 - \theta^*_2)(\theta^*_1 - \theta^*_2)} \\
\end{array}
\right).
$$

\item The transition matrix from the basis (\ref{eq:U2bases}) to the basis (\ref{eq:Astarbases}) is \\
$$ \left(
\begin{array}{ c c c c }
1 &  \frac{\phi}{(\theta^*_1 - \theta^*_0)(\theta^*_1 - \theta^*_2)} & \frac{\varphi + \varphi_2(\theta_0 - \theta_2)(\theta^*_0 - \theta^*_1)}{(\theta^*_1 - \theta^*_0)(\theta^*_1 - \theta^*_2)} & \phi\\
0 & 0 & \theta_2 - \theta_0 & 0  \\
0  & (\theta^*_1 - \theta^*_2)^{-1}  &  (\theta^*_1 - \theta^*_2)^{-1}   & \theta^*_2 - \theta^*_0  \\
0 & 0  &  0    &  (\theta^*_2 - \theta^*_0)(\theta^*_2 - \theta^*_1) \\
\end{array}
\right),
$$

and the transition matrix from the basis  (\ref{eq:Astarbases}) to the basis (\ref{eq:U2bases}) is \\
$$ \left(
\begin{array}{ c c c c }
1 &  \frac{\phi_1}{\theta^*_0 - \theta^*_1} & \frac{\phi}{\theta^*_0 - \theta^*_1} & \frac{\phi}{(\theta^*_0 - \theta^*_1)(\theta^*_0 - \theta^*_2)}\\
0 & (\theta_0 - \theta_2)^{-1} & \theta^*_1 - \theta^*_2 & 1  \\
0  & (\theta_2 - \theta_0)^{-1}  &   0  & 0  \\
0 & 0  &  0    &  \frac{1}{(\theta^*_0 - \theta^*_2)(\theta^*_1 - \theta^*_2)} \\
\end{array}
\right).
$$

\item The transition matrix from the basis (\ref{eq:U3bases}) to the basis (\ref{eq:Astarbases}) is \\
$$ \left(
\begin{array}{ c c c c }
\phi^{-1} &  \frac{1}{(\theta^*_1 - \theta^*_0)(\theta^*_1 - \theta^*_2)} &\frac{\varphi + \varphi_2(\theta_0 - \theta_2)(\theta^*_0 - \theta^*_1)}{(\theta^*_1 - \theta^*_0)(\theta^*_1 - \theta^*_2)\phi} & 1\\
0 & 0 & (\theta_2 - \theta_0)\phi^{-1} & 0  \\
(\theta^*_0 - \theta^*_2)\varphi^{-1}\phi^{-1}  & \varphi^{-1}(\theta^*_1 - \theta^*_0)^{-1}  &  \phi^{-1}(\theta^*_1 - \theta^*_0)^{-1} & 0  \\
(\theta^*_0 - \theta^*_2)(\theta^*_0 - \theta^*_1)\varphi^{-1}\phi^{-1} & 0 & 0 & 0 \\
\end{array}
\right),
$$

and the transition matrix from the basis  (\ref{eq:Astarbases}) to the basis (\ref{eq:U3bases}) is \\
$$ \left(
\begin{array}{ c c c c }
0 &  0 & 0 & \frac{\varphi \phi}{(\theta^*_0 - \theta^*_1)(\theta^*_0 - \theta^*_2)} \\
0 & \varphi(\theta_0 - \theta_2)^{-1} & (\theta^*_1 - \theta^*_0)\varphi & \varphi  \\
0  & \phi(\theta_2 - \theta_0)^{-1}  &  0 & 0  \\
1 & \frac{\varphi_2}{\theta^*_2 - \theta^*_1} &\frac{\varphi}{\theta^*_2 - \theta^*_1}& \frac{\varphi}{(\theta^*_2 - \theta^*_1)(\theta^*_2 - \theta^*_0)} \\
\end{array}
\right).
$$

\item The transition matrix from the basis (\ref{eq:U4bases}) to the basis (\ref{eq:Astarbases}) is \\
$$ \left(
\begin{array}{ c c c c }
\varphi^{-1} &\frac{\varphi + \varphi_1(\theta_0 - \theta_2)(\theta^*_1 - \theta^*_2)}{(\theta^*_1 - \theta^*_0)(\theta^*_1 - \theta^*_2)\varphi} &  \frac{1}{(\theta^*_1 - \theta^*_0)(\theta^*_1 - \theta^*_2)} & 1\\
0 &  (\theta_0 - \theta_2)\varphi^{-1} & 0 & 0  \\
(\theta^*_0 - \theta^*_2)\varphi^{-1}\phi^{-1}  & \varphi^{-1}(\theta^*_1 - \theta^*_0)^{-1}  &  \phi^{-1}(\theta^*_1 - \theta^*_0)^{-1} & 0  \\
(\theta^*_0 - \theta^*_2)(\theta^*_0 - \theta^*_1)\varphi^{-1}\phi^{-1} & 0 & 0 & 0 \\
\end{array}
\right),
$$

and the transition matrix from the basis  (\ref{eq:Astarbases}) to the basis (\ref{eq:U4bases}) is \\
$$ \left(
\begin{array}{ c c c c }
0 &  0 & 0 & \frac{\varphi \phi}{(\theta^*_1 - \theta^*_0)(\theta^*_2 - \theta^*_0)} \\
0 & \varphi(\theta_0 - \theta_2)^{-1} & 0 & 0  \\
0  & \phi(\theta_2 - \theta_0)^{-1}  &  (\theta^*_1 - \theta^*_0)\phi & \phi \\
1 & \frac{\phi_2}{\theta^*_2 - \theta^*_1} &\frac{\phi}{\theta^*_2 - \theta^*_1}& \frac{\phi}{(\theta^*_2 - \theta^*_1)(\theta^*_2 - \theta^*_0)} \\
\end{array}
\right).
$$
\end{enumerate}

\noindent \rm {\bf Proof}. (i) We start with the transition matrix from (\ref{eq:Ubases}) to (\ref{eq:Astarbases}). To get the first column of this matrix, note that the first basis vectors in (\ref{eq:Ubases}) and (\ref{eq:Astarbases}) are the same.\\

\noindent To obtain the second column we observe the
following. 
\noindent By row $[0^*D]$ of Lemma \ref{lem:actionsplit},
\begin{eqnarray}
\label{eq:0starDstareta0star}
(A^* - \theta^*_0I)\eta^*_0 = 0.
\end{eqnarray}

\noindent By Lemma \ref{lem:E0stareigen}(ii),
\begin{eqnarray}
\label{eq:longphieta0star}
(A^* - \theta^*_1I)(A^* - \theta^*_2I)(A - \theta_1I)(A - \theta_2I)\eta^*_0 = \phi \eta^*_0.
\end{eqnarray}

\noindent By Lemma \ref{lem:releigenseq}(i),
\begin{eqnarray}
\label{eq:longvarphi1eta0star}
(A^* - \theta^*_1I)(A - \theta_0I)\eta^*_0 = \varphi_1 \eta^*_0.
\end{eqnarray}

\noindent By Definition \ref{def:initvect},
\begin{eqnarray}
\label{eq:eta0startoeta2}
(A - \theta_1I)(A - \theta_0I)\eta^*_0 = \eta_2.
\end{eqnarray}

\noindent  Consider the equation which is $(\theta^*_1 - \theta^*_0)^{-1}(\theta^*_1 - \theta^*_2)^{-1}$ times (\ref{eq:longphieta0star}) plus $(\theta_0 - \theta_2)(\theta^*_1 - \theta^*_2)^{-1}$ times (\ref{eq:longvarphi1eta0star}).
Adding 
\begin{eqnarray*}
(\theta^*_1 - \theta^*_2)^{-1}(A^* - \theta^*_2I)\eta_2  + (\theta_0 - \theta_2)(A - \theta_0I)\eta^*_0 + (\theta_0 - \theta_2)(\theta_0 - \theta_1)(\theta^*_0 - \theta^*_2)(\theta^*_1 - \theta^*_2)^{-1}\eta^*_0 
\end{eqnarray*}
to both sides of this equation and simplifying the result using (\ref{eq:0starDstareta0star}), (\ref{eq:eta0startoeta2}) and the expression for $E^*_1$ in Notation \ref{note:notephi} we routinely obtain
\begin{eqnarray*}
E^*_1\eta_0 = \frac{\phi +\phi_2(\theta_0 - \theta_2)(\theta^*_1 - \theta^*_0)}{(\theta^*_1 - \theta^*_0)(\theta^*_1 - \theta^*_2)}\eta^*_0 + (\theta_0 - \theta_2)(A-\theta_0I)\eta^*_0 + \frac{(A^* - \theta^*_2I)\eta_2}{\theta^*_1 - \theta^*_2}.
\end{eqnarray*}

\medskip
\noindent To obtain the third column, we observe by Lemma \ref{lem:E0stareigen}(i) and Definition \ref{def:initvect} that
\begin{eqnarray*}
E^*_1\eta_2 &=& \frac{(A^* - \theta^*_0I)(A^* - \theta^*_2I)}{(\theta^*_1 - \theta^*_0)(\theta^*_1 - \theta^*_2)}\eta_2\\
&=&\frac{(A^* - \theta^*_1I)(A^* - \theta^*_2I)(A-\theta_1I)(A-\theta_0I)\eta^*_0 + (\theta^*_1-\theta^*_0)(A^* - \theta^*_2I)\eta_2}{(\theta^*_1 - \theta^*_0)(\theta^*_1 - \theta^*_2)}\\
&=& \frac{\varphi}{(\theta^*_1 - \theta^*_0)(\theta^*_1 - \theta^*_2)}\eta^*_0 + \frac{(A^* - \theta^*_2I)\eta_2}{\theta^*_1 - \theta^*_2}.
\end{eqnarray*}

\medskip
\noindent To obtain the fourth column we use Lemma \ref{lem:E0stareigen}(i) and Definition \ref{def:initvect} to get
\begin{eqnarray*}
\eta^*_2 &=& (A^* - \theta^*_1I)(A^* - \theta^*_0I)\eta_2\\
&=& \varphi\eta^*_0 + (\theta^*_2 - \theta^*_0)(A^* - \theta^*_1I)(A-\theta_1I)(A-\theta_0I)\eta^*_0\\
&=& \varphi\eta^*_0 + (\theta^*_2 - \theta^*_0)(A^* - \theta^*_2I)\eta_2 + (\theta^*_2 - \theta^*_0)(\theta^*_2 - \theta^*_1)\eta_2.
\end{eqnarray*}

\noindent We have now obtained the transition matrix from (\ref{eq:Ubases}) to (\ref{eq:Astarbases}).\\\\
\noindent The transition matrix from the basis  (\ref{eq:Astarbases}) to the basis (\ref{eq:Ubases}) is as shown, since it is routine to verify that its product with the previous matrix  is the identity.\\\\
\noindent (ii) To obtain the transition matrix from the basis (\ref{eq:U2bases}) to the basis (\ref{eq:Astarbases}), compute the product of the transition matrix from the basis (\ref{eq:U2bases}) to the basis (\ref{eq:Ubases}) given in Theorem \ref{thm:TransMat}(i) and the transition matrix from the basis (\ref{eq:Ubases}) to the basis (\ref{eq:Astarbases}) given in Theorem \ref{thm:TransMat4}(i). Simplify the product using (\ref{eq:lambda})--(\ref{eq:mu}). To obtain the transition matrix from the basis (\ref{eq:Astarbases}) to the basis (\ref{eq:U2bases}), compute the product of the transition matrix from the basis (\ref{eq:Astarbases}) to the basis (\ref{eq:Ubases}) given in Theorem \ref{thm:TransMat4}(i) and the transition matrix from the basis (\ref{eq:Ubases}) to the basis (\ref{eq:U2bases}) given in Theorem \ref{thm:TransMat}(i). Simplify the product using (\ref{eq:lambda})--(\ref{eq:mu}).\\\\
\noindent (iii) To obtain the transition matrix from the basis (\ref{eq:U3bases}) to the basis (\ref{eq:Astarbases}), compute the product of the transition matrix from the basis (\ref{eq:U3bases}) to the basis (\ref{eq:Ubases}) given in Theorem \ref{thm:TransMat2}(i) and the transition matrix from the basis (\ref{eq:Ubases}) to the basis (\ref{eq:Astarbases}) given in Theorem \ref{thm:TransMat4}(i). Simplify the product using (\ref{eq:lambda})--(\ref{eq:mu}). To obtain the transition matrix from the basis (\ref{eq:Astarbases}) to the basis (\ref{eq:U3bases}), compute the product of the transition matrix from the basis (\ref{eq:Astarbases}) to the basis (\ref{eq:Ubases}) given in Theorem \ref{thm:TransMat4}(i) and the transition matrix from the basis (\ref{eq:Ubases}) to the basis (\ref{eq:U3bases}) given in Theorem \ref{thm:TransMat2}(i). Simplify the product using (\ref{eq:lambda})--(\ref{eq:mu}).\\\\
\noindent (iv)  To obtain the transition matrix from the basis (\ref{eq:U4bases}) to the basis (\ref{eq:Astarbases}), compute the product of the transition matrix from the basis (\ref{eq:U4bases}) to the basis (\ref{eq:Ubases}) given in Theorem \ref{thm:TransMat}(iv) and the transition matrix from the basis (\ref{eq:Ubases}) to the basis (\ref{eq:Astarbases}) given in Theorem \ref{thm:TransMat4}(i). Simplify the product using (\ref{eq:lambda})--(\ref{eq:mu}). To obtain the transition matrix from the basis (\ref{eq:Astarbases}) to the basis (\ref{eq:U4bases}), compute the product of the transition matrix from the basis (\ref{eq:Astarbases}) to the basis (\ref{eq:Ubases}) given in Theorem \ref{thm:TransMat4}(i) and the transition matrix from the basis (\ref{eq:Ubases}) to the basis (\ref{eq:U4bases}) given in Theorem \ref{thm:TransMat}(iv). Simplify the product using (\ref{eq:lambda})--(\ref{eq:mu}). \hfill $\Box$
\end{theorem}

\bigskip
\noindent We now display the transition matrices between our eigenbasis for $A$ and our eigenbasis for $A^*$.
\medskip
\begin{theorem} \label {thm:TransMat5}
With reference to Notation \ref{note:notephi} the transition matrix from the basis (\ref{eq:Abases}) to the basis (\ref{eq:Astarbases}) is \\
{\small $$ \left(
\begin{array}{ c c c c }
\frac{1}{(\theta_0 - \theta_1)(\theta_0 - \theta_2)} &  \frac{\phi + \phi_2(\theta_0 - \theta_2)(\theta^*_1 - \theta^*_0)}{(\theta_0 - \theta_1)(\theta_0 - \theta_2)(\theta^*_1 - \theta^*_0)(\theta^*_1 - \theta^*_2)} & \frac{\varphi}{(\theta_0 - \theta_1)(\theta_0 - \theta_2)(\theta^*_1 - \theta^*_0)(\theta^*_1 - \theta^*_2)} & \frac{\varphi}{(\theta_0 - \theta_1)(\theta_0 - \theta_2)} \\
1 & \frac{\phi}{(\theta^*_0 - \theta^*_2)(\theta^*_1 - \theta^*_0)} & \frac{\varphi}{(\theta^*_0 - \theta^*_2)(\theta^*_1 - \theta^*_0)} & 0\\
0 & \frac{1}{(\theta^*_2 - \theta^*_1)(\theta^*_0 - \theta^*_2)} & \frac{1}{(\theta^*_2 - \theta^*_1)(\theta^*_0 - \theta^*_2)} & 1\\
\frac{1}{(\theta_0 - \theta_2)(\theta_1 - \theta_2)} &  \frac{\phi}{(\theta_0 - \theta_2)(\theta_1 - \theta_2)(\theta^*_1 - \theta^*_0)(\theta^*_1 - \theta^*_2)}  & \frac{\varphi + \varphi_2(\theta_0 - \theta_2)(\theta^*_0 - \theta^*_1)}{(\theta_0 - \theta_2)(\theta_1 - \theta_2)(\theta^*_1 - \theta^*_0)(\theta^*_1 - \theta^*_2)} & \frac{\phi}{(\theta_0 - \theta_2)(\theta_1 - \theta_2)}   \\
\end{array}
\right),
$$}

and the transition matrix from the basis  (\ref{eq:Astarbases}) to the basis (\ref{eq:Abases}) is \\
{\small $$ \left(
\begin{array}{ c c c c }
\frac{\phi}{(\theta^*_0 - \theta^*_1)(\theta^*_0 - \theta^*_2)} &  \frac{\varphi + \varphi_1(\theta_1 - \theta_2)(\theta^*_0 - \theta^*_2)}{(\theta_1 - \theta_0)(\theta_1 - \theta_2)(\theta^*_0 - \theta^*_1)(\theta^*_0 - \theta^*_2)} & \frac{\varphi\phi}{(\theta_1 - \theta_0)(\theta_1 - \theta_2)(\theta^*_0 - \theta^*_1)(\theta^*_0 - \theta^*_2)} & \frac{\varphi}{(\theta^*_0 - \theta^*_1)(\theta^*_0 - \theta^*_2)} \\
1 & \frac{1}{(\theta_1 - \theta_0)(\theta_0 - \theta_2)} & \frac{\varphi}{(\theta_1 - \theta_0)(\theta_0 - \theta_2)} & 0\\
0 &  \frac{1}{(\theta_2 - \theta_1)(\theta_0 - \theta_2)} & \frac{\phi }{(\theta_2 - \theta_1)(\theta_0 - \theta_2)} & 1 \\\frac{1}{(\theta^*_1 - \theta^*_2)(\theta^*_0 - \theta^*_2)} & \frac{1}{(\theta_1 - \theta_0)(\theta_1 - \theta_2)(\theta^*_0 - \theta^*_2)(\theta^*_1 - \theta^*_2)} & \frac{\varphi + \varphi_2(\theta_1 - \theta_0)(\theta^*_2 - \theta^*_0)}{(\theta_1 - \theta_0)(\theta_1 - \theta_2)(\theta^*_0 - \theta^*_2)(\theta^*_1 - \theta^*_2)} &  \frac{1}{(\theta^*_0 - \theta^*_2)(\theta^*_1 - \theta^*_2)}\\
\end{array}
\right).
$$}

\bigskip
\noindent \rm {\bf Proof}. To obtain the transition matrix from the basis (\ref{eq:Abases}) to the basis (\ref{eq:Astarbases}), compute the product of the transition matrix from (\ref{eq:Abases}) to (\ref{eq:Ubases}) given in Theorem \ref{thm:TransMat3}(i) and the transition matrix from (\ref{eq:Ubases}) to (\ref{eq:Astarbases}) given in Theorem \ref{thm:TransMat4}(i). Simplify the result using (\ref{eq:lambda})--(\ref{eq:mu}). To obtain the transition matrix from (\ref{eq:Astarbases}) to (\ref{eq:Abases}), compute the product of the transition matrix from (\ref{eq:Astarbases}) to (\ref{eq:Ubases}) given in Theorem \ref{thm:TransMat4}(i) and the transition matrix from (\ref{eq:Ubases}) to (\ref{eq:Abases}) given in Theorem \ref{thm:TransMat3}(i). Simplify the result using (\ref{eq:lambda})--(\ref{eq:mu}).  \hfill $\Box$\\
\end{theorem}

\medskip
\noindent We are now ready to prove Theorem \ref{thm:AAstarbasis}.\\

\medskip
\noindent {\bf Proof of Theorem \ref{thm:AAstarbasis}}
\noindent (i) We use (\ref{eq:lambda})--(\ref{eq:mu}) to routinely verify that the given matrices are $T^{-1}BT$ and $T^{-1}B^*T$, where $B$ (respectively $B^*$) denotes the matrix representing $A$ (respectively $A^*$) with respect to the basis (\ref{eq:Ubases}), and $T$ denotes the transition matrix from the basis (\ref{eq:Ubases}) to the basis (\ref{eq:Abases}).\\
\noindent (ii) We use (\ref{eq:lambda})--(\ref{eq:mu}) to routinely verify that the given matrices are $T^{-1}BT$ and $T^{-1}B^*T$, where $B$ (respectively $B^*$) denotes the matrix representing $A$ (respectively $A^*$) with respect to the basis (\ref{eq:Ubases}), and $T$ denotes the transition matrix from the basis (\ref{eq:Ubases}) to the basis (\ref{eq:Astarbases}).\\

\medskip
\section{The classification of TD pairs of shape $(1, 2, 1)$}
\begin{theorem} \label{TDclassif} Given a sequence of scalars
\begin{eqnarray}
\label{TDSysPA}(\lbrace{\theta_i\rbrace}_{i=0}^2, \lbrace{\theta^*_i\rbrace}_{i=0}^2, \varphi, \phi)
\end{eqnarray}
\noindent taken from $\mathbb F$, there exists a TD system $\Phi$ over $\mathbb F$ of shape $(1, 2, 1)$ and parameter array (\ref{TDSysPA}) if and only if (i)--(iii) hold below:

\begin{enumerate}
\item $\theta_i \neq \theta_j$, \qquad $\theta^*_i \neq \theta^*_j$ \qquad if $i \neq j$ \qquad $(0 \leq i, j \leq 2)$;
\item $\varphi \neq 0$, \qquad $\phi \neq 0$;
\item $\varphi \neq \varphi_1 \varphi_2$, \quad where
\begin{eqnarray*}
\varphi_1:= \frac{\phi - \varphi}{(\theta_0 - \theta_2)(\theta^*_0 - \theta^*_2)} -  (\theta_0 - \theta_1)(\theta^*_0 - \theta^*_1),\\
\varphi_2:= \frac{\phi - \varphi}{(\theta_0 - \theta_2)(\theta^*_0 - \theta^*_2)} -  (\theta_1 - \theta_2)(\theta^*_1 - \theta^*_2).
\end{eqnarray*}
\end{enumerate}

\noindent \rm {\bf Proof}. To prove the theorem in one direction let $\Phi:= (A; {\lbrace{E_i}\rbrace}_{i=0}^2; A^*;{\lbrace{E_i^*}\rbrace}_{i=0}^2)$ denote a TD system over $\mathbb F$ of shape $(1, 2, 1)$ and parameter array (\ref{TDSysPA}). We verify conditions (i)--(iii) of the theorem. Condition (i) follows from Definition \ref{def:eigenseq} and condition (ii) follows from Lemma \ref{lem:E0stareigen}. To verify (iii) assume $\varphi = \varphi_1 \varphi_2$. We display a one-dimensional subspace $W$ of $V$ such that $AW \subseteq W$ and $A^*W \subseteq W$. Consider the vector

\begin{eqnarray*}
w = \varphi_2(A - \theta_0I)\eta^*_0 - (A^* - \theta^*_2I)\eta_2,
\end{eqnarray*}

\noindent where $\eta^*_0$ and $\eta_2$ are from Definition \ref{def:initvect}. Note that $w \neq 0$ by Lemma \ref{lem:Ubases}(iv). Using $\varphi = \varphi_1 \varphi_2$ and Theorem \ref{thm:Asplit}(i) we find $Aw = \theta_1w$ and $A^*w = \theta^*_1w$. By these comments $W = \mbox {\rm Span}(w)$ is a one-dimensional subspace of $V$ such that
$AW \subseteq W$ and $A^*W \subseteq W$. This contradicts Definition \ref{def:TDsys}(vi) so $\varphi \neq \varphi_1 \varphi_2$. We have now verified condition (iii) and the theorem is proved in one direction.\\\\
\noindent To prove the theorem in the other direction, we assume the scalars (\ref{TDSysPA}) satisfy (i)--(iii) and display a TD system $\Phi$ over $\mathbb F$ that has shape $(1, 2, 1)$ and parameter array (\ref{TDSysPA}). Let $V$ denote the vector space $\mathbb F^4$ (column vectors). Define the matrices\\

$$ A = \left(
\begin{array}{ c c c c }
\theta_0 & 0  &  0    & 0  \\
1 & \theta_1  &  0    & 0  \\
0 & 0  &  \theta_1 &    0  \\
0 & 1  &  \varphi_2    & \theta_2  \\
\end{array}
\right), \qquad
A^* = \left(
\begin{array}{ c c c c }
\theta^*_0 & \varphi_1  &  \varphi    & 0  \\
0 & \theta^*_1  &  0   &  0    \\
0  & 0  & \theta^*_1   & 1 \\
0  & 0  & 0  & \theta^*_2 \\
\end{array}
\right)$$\\
\noindent and view $A$, $A^*$ as linear transformations on $V$. We show $A$ is diagonalizable. Let $T$ denote the first matrix in Theorem \ref{thm:TransMat3}(i). One checks that $AT = TD$ where $D = \mbox{\rm diag}(\theta_0, \theta_1, \theta_1, \theta_2)$. One also checks that $T$ is invertible; therefore $T^{-1}AT = D$ so $A$ is diagonalizable. We show $A^*$ is diagonalizable. Let $S$ denote the first matrix in Theorem \ref{thm:TransMat4}(i). One checks that $A^*S = SD^*$ where $D^* = \mbox{\rm diag}(\theta^*_0, \theta^*_1, \theta^*_1, \theta^*_2)$. One also checks that $S$ is invertible; therefore $S^{-1}A^*S = D^*$ so $A^*$ is diagonalizable. From the construction the scalars $\theta_0$, $\theta_1$, $\theta_2$ (respectively $\theta^*_0$, $\theta^*_1$, $\theta^*_2$) are  the eigenvalues of $A$ (respectively $A^*$); let $E_0, E_1, E_2$ (respectively $E^*_0, E^*_1, E^*_2$) denote the corresponding primitive idempotents. We show $\Phi = (A; {\lbrace{E_i}\rbrace}_{i=0}^2; A^*;{\lbrace{E_i^*}\rbrace}_{i=0}^2)$ is a TD system over $\mathbb F$ that has shape $(1, 2, 1)$ and parameter array (\ref{TDSysPA}). To verify that $\Phi$ is a TD system, we show that $\Phi$ satisfies the conditions (i)--(vi) of Definition \ref{def:TDsys}. We already verified condition (i) and conditions (ii), (iii) hold by the construction. Next we verify condition (iv). Since $d = 2$, we only need to show that $E_0A^*E_2 = E_2A^*E_0 = 0$. Setting $d = 2$ and $i = 0$, $i = 2$ in
(\ref{eq:Eidef}),\\
$$ E_0 = \left(
\begin{array}{ c c c c }
1 & 0  &  0    & 0  \\
\frac{1}{\theta_0 - \theta_1} & 0  &  0    & 0  \\
0 & 0  &  0 &    0  \\
\frac{1}{(\theta_0 - \theta_1)(\theta_0 - \theta_2)} & 0  &  0    & 0  \\
\end{array}
\right), \qquad
E_2 = \left(
\begin{array}{ c c c c }
0 & 0  &  0    & 0  \\
0 & 0  &  0   &  0    \\
0  & 0  & 0   & 0 \\
\frac{1}{(\theta_2 - \theta_0)(\theta_2 - \theta_1)}  & \frac{1}{\theta_2 - \theta_1}  & \frac{\varphi_2}{\theta_2 - \theta_1}  & 1 \\
\end{array}
\right).$$\\
\noindent Using this data it is routine to verify that $E_0A^*E_2 = E_2A^*E_0 = 0$. We have now verified condition (iv). To verify condition (v) we similarly show $E^*_0AE^*_2 = E^*_2AE^*_0 = 0$. By (\ref{eq:Eidef}),\\
$$E^*_0 = \left(
\begin{array}{ c c c c }
1 & \frac{\varphi_1}{\theta^*_0 - \theta^*_1}  &  \frac{\varphi}{\theta^*_0 - \theta^*_1}    & \frac{\varphi}{(\theta^*_0 - \theta^*_1)(\theta^*_0 - \theta^*_2)}  \\
0 & 0  &  0  & 0  \\
0 & 0  &  0  & 0  \\
0 & 0  &  0  & 0  \\
\end{array}
\right), \qquad
E^*_2 = \left(
\begin{array}{ c c c c }
0 & 0  &  0    & \frac{\varphi}{(\theta^*_2 - \theta^*_0)(\theta^*_2 - \theta^*_1)}  \\
0 & 0  &  0   &  0    \\
0  & 0  & 0   & \frac{1}{\theta^*_2 - \theta^*_1} \\
0  & 0  & 0  & 1 \\
\end{array}
\right).$$\\
\noindent By the above line it is routine to verify that $E^*_0AE^*_2 = E^*_2AE^*_0 = 0$. We have now verified condition (v). Finally, we claim that condition (vi) holds. To prove the claim we start with some comments. By construction the span of the vector $(1, 0, 0, 0)^t$ is $E^*_0V$, the span of $(0, 0, 0, 1)^t$ is $E_2V$, the span of $(0, 1, 0, 0)^t$ is $(A - \theta_0I)E^*_0V$, and the span of $(0, 0, 1, 0)^t$ is $(A^* - \theta^*_2I)E_2V$. Therefore
\begin{eqnarray}
\label{eq:e1}
V= E^*_0V + (A-\theta_0I)E^*_0V + (A^*-\theta^*_2I)E_2V + E_2V  \qquad \qquad (\mbox{\rm direct sum}).
\end{eqnarray}
By the form of $A$,
\begin{eqnarray}
\label{eq:e2}
E_2V= 
(A-\theta_1 I)(A-\theta_0 I)E^*_0V.
\end{eqnarray}
By the form of $A^*$ and since $\varphi\neq 0$,
\begin{eqnarray}
\label{eq:e3}
E^*_0V= 
(A^*-\theta^*_1 I)(A^*-\theta^*_2 I)E_2V.
\end{eqnarray}
Similarly by the form of $E_0, E^*_2$ and since $\phi \neq 0$,
\begin{eqnarray}
\label{eq:e4}
V&=& E^*_0V
+(A-\theta_2I)E^*_0V
+(A^*-\theta^*_2I)E_0V
+E_0V  \quad \quad (\mbox{\rm direct sum}),
\\
\label{eq:e5}
E_0V&=& 
(A-\theta_1 I)(A-\theta_2 I)E^*_0V,
\\
\label{eq:e6}
E^*_0V&=& 
(A^*-\theta^*_1 I)(A^*-\theta^*_2 I)E_0V,
\end{eqnarray}
and
\begin{eqnarray}
\label{eq:e7}
V&=& E^*_2V
+(A-\theta_0I)E^*_2V
+(A^*-\theta^*_0I)E_2V
+E_2V  \quad \quad (\mbox{\rm direct sum}),
\\
\label{eq:e8}
E_2V&=& 
(A-\theta_1 I)(A-\theta_0 I)E^*_2V,
\\
\label{eq:e9}
E^*_2V&=& 
(A^*-\theta^*_1 I)(A^*-\theta^*_0 I)E_2V.
\end{eqnarray}
Let $W$ denote a proper subspace of $V$ that is invariant
under $A$ and $A^*$. We show $W = 0$. Note that $E^*_0W = 0$;
otherwise $E^*_0V\subseteq W$ and then $W = V$ by
(\ref{eq:e1})--(\ref{eq:e3}).
Similarly $E_2W = 0$; otherwise $E_2V\subseteq W$ and then $W = V$ by
(\ref{eq:e1})--(\ref{eq:e3}).
We have $E_0W = 0$; otherwise $E_0V\subseteq W$ and then
$W = V$ by 
(\ref{eq:e4})--(\ref{eq:e6}).
Also $E^*_2W = 0$; otherwise $E^*_2V\subseteq W$ and then
$W = V$ by 
(\ref{eq:e7})--(\ref{eq:e9}).
For $w \in W$ we show $w = 0$. Write $w=(a,b,c,d)^t$.
By our above comments $E_0w = 0$, so $a = 0$ by the form of $E_0$.
Similarly $E^*_2w = 0$, so $d = 0$ by the form of $E^*_2$.
We also have $E^*_0w = 0$, $E_2w = 0$;
evaluating these using $a = 0$, $d = 0$ and the forms of $E^*_0, E_2$
we get
\begin{eqnarray*}
   \varphi_1 b + \varphi c &=&0,\\
   b +   \varphi_2 c &=&0.
\end{eqnarray*}
Solving the above system of equations using
$\varphi \neq \varphi_1 \varphi_2$ we routinely find
$b = 0$ and $c = 0$. Now each of $a, b, c, d$ is $0$ so $w = 0$.
We have now shown $W = 0$ and this establishes condition (vi).
We have shown $\Phi$ satisfies conditions (i)--(vi) of Definition \ref{def:TDsys} so $\Phi$ is a TD system.
By construction $\Phi$ is over $\mathbb F$ and has shape $(1, 2, 1)$. Using Lemma \ref{lem:E0stareigen} and Definition \ref{def:pararray} we routinely find that the sequence (\ref{TDSysPA}) is the parameter array of $\Phi$. \hfill $\Box$
\end{theorem}

\bigskip
\noindent{\Large\bf Acknowledgements}

\medskip
\noindent
This paper was written while the author was an honorary
fellow at the University of Wisconsin-Madison, November 2006--October 2007, supported by the HEDP-FDP Sandwich Program of the Commision on Higher Education, Philippines.
The author would like to thank his UW-Madison advisor Paul Terwilliger
for his many valuable ideas and suggestions. The author is also greatly indebted to his advisor in the Philippines, Arlene Pascasio, who facilitated his application to the above funding agency, and along with Junie Go, gave this paper a careful reading.

\medskip
\noindent
Melvin A. Vidar\\
Math and Statistics Department\\
College of Arts and Sciences\\
University of the East\\
Recto Ave. Manila, Philippines\\
email: melvs\_edu@yahoo.com

\end{document}